\documentclass[reqno, table]{amsart} 
\usepackage{amsmath}
\usepackage[margin=1.0in]{geometry}
\usepackage{amscd,amsmath,amsxtra,amsthm,amssymb,stmaryrd,xr,mathrsfs,mathtools, enumerate,commath, comment, enumitem, graphicx}
\usepackage{amsfonts, amssymb, amsthm, amsmath, braket, hyperref, mathtools, comment, mathrsfs, enumitem, cleveref}
\usepackage{stmaryrd}
\usepackage{amsfonts}
\usepackage{orcidlink}
\usepackage{multirow}
\usepackage{xcolor}
\usepackage{commath}
\usepackage{enumitem}
\usepackage{comment}
\usepackage{tikz-cd}
\usepackage{tkz-graph}
\usepackage{colonequals}
\usepackage{hyperref}
\usepackage{graphicx}
\usepackage{bigstrut}
\usepackage{enumerate}
\usepackage{tabularx}
\numberwithin{equation}{section}
\usepackage{longtable, multirow, array}
\usepackage{cleveref}
\crefname{section}{§}{§§}
\Crefname{section}{§}{§§}
\newcommand{\genlegendre}[4]{%
	\genfrac{(}{)}{}{#1}{#3}{#4}%
	\if\relax\detokenize{#2}\relax\else_{\!#2}\fi
}

\usepackage{longtable} % for 'longtable' environment
\usepackage{pdflscape} % for 'landscape' environment
\usepackage{booktabs}
\usepackage{hyperref}
\definecolor{vegasgold}{rgb}{0.77, 0.7, 0.35}
\definecolor{darkgoldenrod}{rgb}{0.72, 0.53, 0.04}
\definecolor{gold(metallic)}{rgb}{0.83, 0.69, 0.22}
\hypersetup{
 colorlinks=true,
 linkcolor=darkgoldenrod,
 filecolor=brown,      
 urlcolor=gold(metallic),
 citecolor=darkgoldenrod,
 }
\newtheorem{lthm}{Theorem}

\usepackage[all,cmtip]{xy}
\DeclareFontFamily{U}{wncy}{}
\DeclareFontShape{U}{wncy}{m}{n}{<->wncyr10}{}
\DeclareSymbolFont{mcy}{U}{wncy}{m}{n}
\DeclareMathSymbol{\Sh}{\mathord}{mcy}{"58}
\usepackage[T2A,T1]{fontenc}
\usepackage[OT2,T1]{fontenc}

\usepackage{tikz}
\usetikzlibrary{shapes.geometric}
\usetikzlibrary{decorations.markings}
\tikzset{every loop/.style={min distance=10mm,looseness=10}}
\tikzstyle{vertex}=[auto=left,circle,minimum size=1pt,inner sep=0pt]

\newtheorem{theorem}{Theorem}[section]
\newtheorem{lemma}[theorem]{Lemma}

\newtheorem*{theorem*}{Theorem}
\newtheorem*{ass*}{Assumption}

\newtheorem{corollary}[theorem]{Corollary}
\newtheorem{remark}[theorem]{Remark}

\newtheorem{proposition}[theorem]{Proposition}

\newcommand{\Z}{\mathbb{Z}}

\newcommand{\p}{\mathfrak{p}}
\newcommand{\Q}{\mathbb{Q}}
\newcommand{\F}{\mathbb{F}}

\DeclareFontEncoding{OT2}{}{}
\DeclareSymbolFont{cyrletters}{OT2}{wncyr}{m}{n}
\DeclareMathSymbol{\Sha}{\mathalpha}{cyrletters}{"58}

\numberwithin{equation}{section}

\begin{document}

\title{Matrix Criteria for Mordell--Weil Ranks of Quartic Twists of $y^2=x^3+x$ with Root Number $-1$}

\author[Debajyoti De]{Debajyoti De}
\address[De]{Department of Mathematics, IIT Madras, India}
\email{debajyotide20@gmail.com}


\begin{abstract}
We study the Mordell--Weil ranks of quartic twists $E_a:y^2=x^3+ax$ over $\Q$ with global root number $-1$. For $D$ a product of $n\ge2$ distinct odd primes, a $2$-isogeny descent yields unconditional rank bounds for $E_D$, $E_{-D}$, and $E_{-2D}$ in terms of matrices over $\F_2$ defined by additive Legendre symbols. Under suitable congruence and matrix-rank conditions, these curves have rank at most one, and rank one assuming the Mordell--Weil parity conjecture. For $E_{-2D}$, an additional explicit square condition establishes rank one unconditionally. Our results extend earlier two-prime criteria to arbitrary numbers of prime factors and further to fourth-power-free coefficients $H=\prod_{i=1}^n\ell_i^{e_i}$, where $e_i\in\{1,3\}$.
\end{abstract}

\maketitle

\section{Introduction}

The study of rational points on elliptic curves has a long history in number theory. Determining the Mordell--Weil rank of an elliptic curve defined over $\Q$ is a central problem in arithmetic geometry. By the Mordell--Weil theorem, the group of rational points $E(\Q)$ is finitely generated and $E(\Q)\cong E(\Q)_{\mathrm{tors}} \oplus\mathbb Z^{r},$
where $r=\operatorname{rank}E(\Q)$ and $E(\Q)_{\mathrm{tors}}$ is the torsion subgroup of $E(\Q)$. Although this theorem establishes finite generation, determining $r$ for a given curve remains difficult. For curves with a rational point of order two, descent through a $2$-isogeny provides an effective way to obtain upper bounds for the rank from explicit arithmetic conditions. In this paper, we study the family
$$E_a:~ y^2=x^3+ax,\quad a\in \Z\setminus \{0\},$$
whose members have $j$-invariant $1728$ and the rational point $(0,0)$ is of order two. This family represents all quartic twists of $E_1:y^2=x^3+x$ over $\Q$, upto $\Q$-isomorphism. Indeed, for any $a,b\in\mathbb Z\setminus\{0\}$ satisfying $b=au^4$ for some $u\in\Q^\times$, the map
$$E_a\longrightarrow E_b,\quad (x,y)\longmapsto(u^2x,u^3y),
$$
defines an isomorphism over $\Q$.  Thus the twist parameter is naturally considered modulo
fourth powers, and we may work with fourth power free
integral coefficients. The subfamily
$E_{-m^2}:y^2=x^3-m^2x$ with $m\geq 1$ consists of the classical
congruent number curves, whose $2$-Selmer ranks have
been studied using Monsky's matrix; see
\cite[Appendix]{HeathBrown1994}.
In the present paper, we focus on coefficients whose
odd prime divisors occur to odd exponents. Since the
coefficient is considered modulo fourth powers,
these exponents may be taken to belong to $\{1,3\}$.

The ranks of several subfamilies of $E_a$ have been investigated using descent and the arithmetic of associated Diophantine equations. Spearman \cite{Spearman,Spearman2} studied the curves
$E_{-p}:y^2=x^3-px$ and $E_{-2p}:y^2=x^3-2px$, while Fujita and Terai
\cite{Fujita-terai} considered $y^2=x^3-p^kx$ for
$k\in\{1,2,3\}$. Related rank questions for coefficients involving one or two varying primes were studied by Daghigh and Didari \cite{Hasaan-Didari}, Kim \cite{kim2015ranks}, and Mina and Bacani
\cite{Mina-Bacani}. Of particular relevance here is the work of Ghosh, Roy and Sharma \cite{Arkabrata-Bidisha-Richa} on $E_{-2pq}: y^2=x^3-2pqx,$
where $p$ and $q$ are distinct odd primes. They obtained explicit rank bounds and assuming the Mordell--Weil parity conjecture, sets of possible ranks determined by the residue classes of $p,q$ modulo $8$ and their Legendre symbols. They also gave sufficient square conditions for the rank to equal one unconditionally.

The Birch and Swinnerton-Dyer (BSD) conjecture relates
the Mordell--Weil rank of an elliptic curve $E/\Q$
to the order of vanishing of its $L$-function at $s=1$.
In particular, it predicts that
$\operatorname{ord}_{s=1}L(E,s)
=\operatorname{rank}E(\Q)$.
By modularity, $L(E,s)$ admits an analytic continuation
to the whole complex plane, and its completed $L$-function
satisfies
$\Lambda(E,s)=\omega(E)\Lambda(E,2-s)$,
where $\omega(E)\in\{\pm1\}$ is the global root number.
Equivalently, $\omega(E)$ is the product of the local root
numbers over all places of $\Q$:
$\omega(E)=\omega_\infty(E)\displaystyle\prod_p \omega_p(E)$,
where $\omega_\infty(E)=-1$ and $p$ runs over the finite primes.
See \cite[Section~1]{ro} for a discussion of local and
global root numbers. Let $a$ be a nonzero fourth-power-free integer with $4\nmid a$, and consider $E_a:y^2=x^3+ax$.
The formula of Birch and Stephens
\cite[Equations~(9)--(10), p.~296]{BS}, gives
\[
\omega(E_a)=\operatorname{sgn}(a)c_2(a)
\prod_{\substack{p\ge3\\p^2\parallel a}}
\left(\frac{-1}{p}\right),
\]
where $c_2(a)=-1$ if
$a\equiv 3,\,5,\,13,\,15 \pmod{16}$, and $c_2(a)=1$ otherwise.
Here $c_2(a)$ is an auxiliary sign in the simplified
global formula, and $p^2\parallel a$ means that
$p^2$ divides $a$, but $p^3$ does not.

If every odd prime divisor of $a$ occurs with exponent
$1$ or $3$, the product is empty, so
$\omega(E_a)=\operatorname{sgn}(a)c_2(a)$.
This applies to all coefficients occurring in
Theorems~\ref{Thm Positive D}, \ref{Thm Negative D},
\ref{Thm negative 2D}, and
\ref{thm odd prime exponent extension}.
Table~\ref{tab:root-numbers} records the square-free specialization,
with the signed coefficient denoted by $a$.

\begin{table}[htbp]
\centering
\renewcommand{\arraystretch}{1.15}
\begin{tabular}{cclc}
\hline
Sign of $a$ & Parity of $a$
& Congruence condition & $\omega(E_a)$\\
\hline
$a>0$ & Odd
& $a\equiv1,7,9,11\pmod{16}$ & $+1$\\
$a>0$ & Odd
& $a\equiv3,5,13,15\pmod{16}$ & $-1$\\
$a>0$ & Even
& All & $+1$\\
$a<0$ & Odd
& $a\equiv3,5,13,15\pmod{16}$ & $+1$\\
$a<0$ & Odd
& $a\equiv1,7,9,11\pmod{16}$ & $-1$\\
$a<0$ & Even
& All & $-1$\\
\hline
\end{tabular}
\caption{Global root numbers of $E_a:y^2=x^3+ax$
for nonzero square-free integers $a$.}
\label{tab:root-numbers}
\end{table}

In our earlier work \cite{Das-De-Mondal}, we studied rank-zero
questions for the family $E_a:y^2=x^3+ax$ in cases with global
root number $+1$. Under suitable congruence conditions on $a$
and its prime divisors, we used descent through a rational
$2$-isogeny to obtain sufficient conditions for
$\operatorname{rank}E_a(\Q)=0$ in terms of the ranks
of matrices encoding quadratic residue data. These results
cover negative odd, positive odd, and positive even
coefficients. In the present paper, we continue this study for quartic twists with global root number $-1$. The functional equation then forces the analytic rank to be odd. The Mordell--Weil parity conjecture asserts that
\begin{equation*}\label{eq:MW-parity-introduction}
(-1)^{\operatorname{rank}E(\Q)}=\omega(E).
\end{equation*}
Thus, when the global root number is $-1$, an unconditional upper bound of one determines the rank under this parity conjecture. When $\omega(E)=-1$, an unconditional upper bound
$\operatorname{rank}E(\Q)\leq1$ gives rank one
under the Mordell--Weil parity conjecture. Together
with this upper bound, a point of infinite order
proves rank one unconditionally.

Several arithmetic criteria for determining the Mordell--Weil rank of $E_a$ have been developed using descent, Selmer groups, and class groups. Viewing $E_a\colon y^2=x^3+ax$ as a quartic twist of $E_1\colon y^2=x^3+x$, De, Majumdar, and Mondal show that, for each odd prime $p\nmid a$, the $p$-Selmer rank of $E_a$ is determined by its global root number and a suitable character component of the relative $p$-class group of a cyclic quartic extension. The base field of this extension contains $\mathbb{Q}(i)$ and for fixed $p$,
is independent of $a$ up to isomorphism \cite[Theorems~3.1 and~3.2]{DDM}. In the case $\omega(E_a)=-1$ considered here, their results provide a computational criterion for establishing Mordell--Weil rank one
and finiteness of the Tate--Shafarevich group. It suffices to find a single prime $p\equiv 1\pmod{4}$ with $p\nmid a$ for which the $p$-Selmer rank is one \cite[Corollary~3.4]{DDM}.

In this article, we obtain explicit matrix bounds for the ranks of $E_D$,
$E_{-D}$, and $E_{-2D}$, where $D$ is a product of arbitrarily
many distinct odd primes. In particular, our result for
$E_{-2D}$ extends the two prime results of
\cite{Arkabrata-Bidisha-Richa} to arbitrary prime support. The
same framework also gives a uniform square condition for
unconditional rank one. Finally, we extend these bounds to $E_H$, $E_{-H}$, and $E_{-2H}$, where $H$ is a positive odd integer with at least two distinct prime factors, each
occurring with exponent $1$ or $3$.

We now introduce some notation that will be used throughout the article.

\textbf{Notation:} We use the following notation. We write $\F_2$ for the field with two elements.
Unless another field is specified, all vector and matrix
operations are over $\F_2$.

\begin{itemize}
\item For a field $K$ and positive integers $m,n$,
$M_{m,n}(K)$ denotes the set of $m\times n$ matrices over $K$,
and $M_n(K)=M_{n,n}(K)$.

\item For $r\geq1$, let
$\mathbf0_r=(0,\ldots,0)^{\mathsf T}$ and
$\mathbf1_r=(1,\ldots,1)^{\mathsf T}$ in $\F_2^r$.

 \item $I_r$ denotes the identity matrix of order $r$ over $\F_2$.
 \item For a vector $\mathbf x=(x_1,\ldots,x_r)^{\mathsf T}$, we denote its $i$-th entry by $(\mathbf x)_i=x_i$, for $1\le i\le r$.

\item For an odd prime $\ell$ and an integer $a$ with
$\ell\nmid a$, define the additive Legendre symbol by
\[
\left[\frac{a}{\ell}\right]=
\begin{cases}
0,&\displaystyle\left(\frac{a}{\ell}\right)=1,\\
1,&\displaystyle\left(\frac{a}{\ell}\right)=-1,
\end{cases}
\qquad\text{in }\F_2,
\]
where $\left(\frac{a}{\ell}\right)$ is the usual
Legendre symbol.

\item For $L=p_1\cdots p_r$, with distinct odd primes
listed in a fixed order, define
$\mathcal A_L=(a_{ij})\in M_r(\F_2)$ by
\begin{equation}\label{matrix A_D}
a_{ij}=\left[\frac{p_j}{p_i}\right]\quad(i\ne j),
\quad
a_{ii}=\sum_{j\ne i}a_{ij}.
\end{equation}

\item Let $L=p_1\cdots p_r$, with distinct odd primes
listed in a fixed order. For a positive divisor
$\delta=\displaystyle\prod_{i=1}^rp_i^{\varepsilon_i(\delta)}$ of $L$, define
\begin{equation}\label{varphi-def}
\varphi(\delta)=
(\varepsilon_1(\delta),\ldots,\varepsilon_r(\delta))^{\mathsf T}
\in\F_2^r.
\end{equation}
Thus $\varphi(1)=\mathbf0_r$ and
$\varphi(L)=\mathbf1_r$.
\end{itemize}

We are now in a position to state our main
results. These results provide explicit upper
bounds for the Mordell--Weil ranks of the
elliptic curves under consideration in terms
of the ranks of matrices over $\F_2$
constructed from quadratic-residue data.

\begin{lthm}\label{Thm Positive D}
Let $D=\ell_1\cdots\ell_n$ be a product of $n\geq2$ distinct
odd primes, and let $\mathcal A_D$ be defined by
\eqref{matrix A_D}. Suppose that
$D\equiv3,5,13,$ or $15\pmod{16}$. Put
$k=n-\operatorname{rank}_{\F_2}\mathcal A_D$.
Then the elliptic curve $E_D:y^2=x^3+Dx$ satisfies
$\omega(E_D)=-1, \text{ and } \operatorname{rank}E_D(\Q)\leq2k-1.$

In particular, if
$\operatorname{rank}_{\F_2}\mathcal A_D=n-1$,
then $\operatorname{rank}E_D(\Q)\leq1$,
with equality if the Mordell--Weil parity conjecture
holds for $E_D$.
\end{lthm}

\begin{lthm}\label{Thm Negative D}
Let $D=\ell_1\cdots\ell_n$ be a product of $n\geq2$ distinct
odd primes, and let $\mathcal A_D$ be defined by
\eqref{matrix A_D}. Suppose that
$D\equiv5,7,9,$ or $15\pmod{16}$. Put
$k=n-\operatorname{rank}_{\F_2}\mathcal A_D$.
Then the elliptic curve $E_{-D}:y^2=x^3-Dx$ satisfies
$\omega(E_{-D})=-1, \text{ and }\operatorname{rank}E_{-D}(\Q)\leq2k-1.$

In particular, if
$\operatorname{rank}_{\F_2}\mathcal A_D=n-1$,
then $\operatorname{rank}E_{-D}(\Q)\leq1$,
with equality if the Mordell--Weil parity conjecture
holds for $E_{-D}$.
\end{lthm}

We now turn to the family $E_{-2D}$, for which we introduce the following augmented matrix.
For $1\leq i\leq n$, put
$b_i=\left[\frac{2}{\ell_i}\right]$ and
$s_i=\left[\frac{-1}{\ell_i}\right]$, and write
$\mathbf s=(s_1,\ldots,s_n)^{\mathsf T}$.
Define
\begin{equation}\label{eq:augmented-RD}
\mathcal R_D=
\left[
\mathcal A_D+\operatorname{diag}(b_1+s_1,\ldots,b_n+s_n)
\,\middle|\,\mathbf s
\right]
\in M_{n,n+1}(\F_2),
\end{equation}
and put
$\rho_D=\operatorname{rank}_{\F_2}\mathcal R_D$.

\begin{lthm}\label{Thm negative 2D}
Let $D=\ell_1\cdots\ell_n$ be a product of $n\geq2$ distinct
odd primes, and let $\mathcal R_D$ and $\rho_D$ be as above.
Set $E_{-2D}: y^2=x^3-2Dx,\text{ and } r_D=\operatorname{rank}E_{-2D}(\Q).$
Then the following statements hold.
\begin{enumerate}
\item The global root number is $\omega(E_{-2D})=-1$, and
the unconditional bound is $r_D\leq2(n-\rho_D)+1$.
In particular, if $\rho_D=n$, then $r_D\leq1$.

\item If the Mordell--Weil parity conjecture holds for
$E_{-2D}$, then $r_D$ is odd and $1\leq r_D\leq2(n-\rho_D)+1$.
In particular, if $\rho_D=n$, then $r_D=1$.

\item Suppose that $\rho_D=n$. If there exist a positive
divisor $\delta\mid D$, a sign $\varepsilon\in\{1,-1\}$,
and positive integers $h,m$ such that $m^2=\varepsilon\left(
\delta-2h^4\frac{D}{\delta}\right),$
then $r_D=1$ without assuming the parity conjecture.
\end{enumerate}
\end{lthm}

When $D=pq$, the matrix $\mathcal R_D$ has two rows and
three columns. Computing its rank recovers the two-prime
rank bounds and conditional rank restrictions in
\cite[Theorems~1.1--1.4]{Arkabrata-Bidisha-Richa} and
\cite[Theorem~1]{MinaBacani2026}.

% \begin{corollary}[The two-prime specialization]
% \label{cor-two-prime-specialization}
% Let $p$ and $q$ be distinct odd primes, and put
% $$E=E_{-2pq}:~ y^2=x^3-2pqx,
% \quad r=\operatorname{rank}E(\Q),
% \quad \chi=\left(\frac pq\right).$$
% Then $\omega(E)=-1$. Assuming the Mordell--Weil parity conjecture
% for $E$, the rank $r$ belongs to the set in the corresponding
% entry of the following table, where residue pairs are listed
% up to interchange:
% \[
% \begin{array}{c|c|c}
% (p,q)\bmod8&\chi=1&\chi=-1\\ \hline
% (1,1)&\{1,3,5\}&\{1,3\}\\[2pt]
% (1,3),(1,5),(1,7)&\{1,3\}&\{1\}\\[2pt]
% \substack{(3,3),(3,5),(3,7),\\(5,5),(5,7)}
% &\{1\}&\{1\}\\[4pt]
% (7,7)&\{1,3\}&\{1,3\}
% \end{array}
% \]
% The largest integer in each displayed set is also an
% unconditional upper bound for $r$.
% \end{corollary}
For $D=pq$, where $p$ and $q$ are distinct
odd primes, the augmented matrix $\mathcal R_D$
has two rows and three columns. Its rank can
therefore be determined explicitly in terms
of the residue classes of $p$ and $q$
modulo $8$ and the Legendre symbol
$\left(\frac{p}{q}\right)$. These calculations recover the conditional
rank restrictions and unconditional upper
bounds obtained in
\cite[Theorems~1.1--1.4]{Arkabrata-Bidisha-Richa}.
The precise two-prime specialization is
stated and proved in
Corollary~\ref{cor-two-prime-specialization}.

Moreover, the square conditions in
\cite[Theorem~1.5]{Arkabrata-Bidisha-Richa} are special cases
of Theorem~\ref{Thm negative 2D}\textup{(3)}. For example,
taking $(\delta,\varepsilon,h)=(pq,1,2)$ gives
$m^2=pq-32$, while taking $(\delta,\varepsilon,h)=(1,-1,1)$
gives $m^2=2pq-1$. The corresponding substitutions for the
remaining residue pairs recover the other criteria. Thus both
the conditional two prime rank restrictions and the stated
unconditional rank-one criteria follow from the same matrix
bound and square condition.

We next allow the odd prime exponents to be $1$ or $3$.
If $H=\displaystyle \prod_i\ell_i^{e_i}$ and $D=\displaystyle\prod_i\ell_i$, then
$H=DF^2$, where $F=\displaystyle\prod_{e_i=3}\ell_i$.
Thus $E_H$ is the quadratic twist of $E_D$ by $F$.
The following theorem gives bounds using the same
prime support matrices, with congruence conditions
imposed on $H$ itself.

\begin{lthm}
\label{thm odd prime exponent extension}
Let $H=\displaystyle\prod_{i=1}^{n}\ell_i^{e_i},~
e_i\in\{1,3\},~ D=\displaystyle\prod_{i=1}^{n}\ell_i,$ where $n\geq2$ and the $\ell_i$ are distinct odd primes.
Define $\mathcal A_D$ and $\mathcal R_D$ by
\eqref{matrix A_D} and \eqref{eq:augmented-RD} respectively, and put $k=n-\operatorname{rank}_{\F_2}\mathcal A_D,
~\rho_D=\operatorname{rank}_{\F_2}\mathcal R_D.$
Then the following statements hold.
\begin{enumerate}
\item If $H\equiv3,\,5,\,13$ or $15\pmod{16}$, then $\omega(E_H)=-1, ~\operatorname{rank}E_H(\Q)\leq2k-1.$

\item If $H\equiv5,\,7,\,9$ or $15\pmod{16}$, then
$\omega(E_{-H})=-1,~ \operatorname{rank}E_{-H}(\Q)\leq2k-1.$

\item For $E_{-2H}:y^2=x^3-2Hx$, one has $\omega(E_{-2H})=-1,~ \operatorname{rank}E_{-2H}(\Q)
\leq2(n-\rho_D)+1.$
\end{enumerate}
Consequently, in parts \textup{(1)} and \textup{(2)}, the
condition $\operatorname{rank}_{\F_2}\mathcal A_D=n-1$
gives rank at most one. In part \textup{(3)}, the condition
$\operatorname{rank}_{\F_2}\mathcal R_D=n$ gives rank
at most one. Under the corresponding matrix hypothesis, if
the Mordell--Weil parity conjecture holds for the curve in
question, its rank equals one.
\end{lthm}

% The proofs use descent through the rational $2$-isogeny
% $E_a\to E_{-4a}$. Quadratic-residue constraints on the
% descent equations yield the matrices $\mathcal A_D$ and
% $\mathcal R_D$. Explicit solutions of these equations
% provide points of infinite order for the unconditional
% rank-one criterion.

\medskip
\noindent\textbf{Organisation of the paper.}
Section~\ref{Preliminaries} establishes the required matrix
identities and local descent conditions and recalls descent
through a rational $2$-isogeny.
Section~\ref{Proof of the main theorem} proves the main results
by encoding necessary quadratic-residue conditions from the
descent equations as linear systems over $\F_2$.
The resulting matrix ranks give unconditional upper bounds
for the Mordell--Weil ranks. These bounds yield rank-one
criteria under the parity conjecture and for $E_{-2D}$,
unconditionally through an explicit square condition.
This section also treats the two-prime specialization and
the extension to odd prime exponents in $\{1,3\}$.
Finally, Section~\ref{Numerical Examples} presents numerical examples.

\subsection*{\it Acknowledgement}  The author acknowledges the Institute Fellowship and the excellent research environment provided by the Indian Institute of Technology Madras. 

\medskip
\noindent\textbf{Declaration of AI assistance}
The author used ChatGPT 6.1 Sol (High) to check grammar,
improve the exposition, and discuss the formulation of
a matrix-rank hypothesis.

Specifically, computations performed by the author using
SageMath initially suggested the sufficient condition
based on the invertibility of
$\mathcal A_D+\operatorname{diag}(b_1,\ldots,b_n)$.
Discussions with ChatGPT 6.1 Sol (High) assisted in
identifying the more general augmented-matrix condition
$\operatorname{rank}_{\F_2}\mathcal R_D=n$
appearing in Theorem~\ref{Thm negative 2D}.

All proof ideas, mathematical arguments, and proofs,
including the proof under the augmented-matrix hypothesis,
were developed independently by the author. The author takes full responsibility for the mathematical correctness,
originality, and final content of the manuscript.

\begin{remark}[Comparison with an invertibility criterion]
Put $A=\mathcal A_D$, $B=\operatorname{diag}(b_1,\ldots,b_n)$,
and $C=A+B$. The matrix identities established below give
\[
\rho_D
=\operatorname{rank}_{\F_2}
   [C^{\mathsf T}+\mathbf s\mathbf s^{\mathsf T}
    \mid\mathbf s]
=\operatorname{rank}_{\F_2}
   [C^{\mathsf T}\mid\mathbf s]
\geq\operatorname{rank}_{\F_2}C.
\]
Consequently, invertibility of $A+B$ implies
$\rho_D=n$ and hence
$\operatorname{rank}E_{-2D}(\Q)\leq1$.
The augmented full-rank hypothesis is strictly
more general.

Indeed, take $D=3\cdot11$, with the primes ordered
as $(3,11)$. Then
$A=\begin{pmatrix}1&1\\0&0\end{pmatrix},~
B=S=\begin{pmatrix}1&0\\0&1\end{pmatrix}.$
Thus $A+B=\begin{pmatrix}0&1\\0&1\end{pmatrix}$
is singular, whereas
$\mathcal R_D
=[A+B+S\mid\mathbf s]
=\begin{pmatrix}1&1&1\\0&0&1\end{pmatrix}$
has rank $2$. Therefore, the augmented-matrix
criterion applies although the invertibility
criterion does not.
\end{remark}

\section{Preliminaries}\label{Preliminaries}

In this section, we recall the basic results on
$2$-descent via $2$-isogeny that will be used
to determine the Mordell--Weil ranks of the
elliptic curves under consideration. We also
establish several linear-algebraic identities
and rank relations over $\F_2$ that
will be needed in the proofs of our main results. We begin with the matrix identities and rank
relations associated with the matrix
$\mathcal A_L$ as defined in \eqref{matrix A_D}.

\subsection{Matrix Identities and Rank Relations}

\begin{lemma}\label{lem:matrix-identities}
Let $L=p_1\cdots p_r$ be a product of distinct odd primes, and
let $A=\mathcal{A}_L$ be the matrix defined in \eqref{matrix A_D}.
For $1\le i\le r$, put $s_i=\bigl[\frac{-1}{p_i}\bigr]$.
Set $\mathbf s=(s_1,\ldots,s_r)^\mathsf T$,
$S=\operatorname{diag}(s_1,\ldots,s_r)$, and $Q=A+S$. Then over $\F_2$ the following identities hold:
$$A \mathbf1_r=\mathbf0_r, \quad A+A^{\mathsf T} =S+\mathbf s\mathbf s^{\mathsf T},~ \text{ and } ~Q=A+S=A^{\mathsf T}+\mathbf s\mathbf s^{\mathsf T}.$$
\end{lemma}

\begin{proof}
All computations are carried out over $\F_2$.
By the definition of the diagonal entries of $A$,
for every $1\le i\le r$, we have
$\displaystyle\sum_{j=1}^r a_{ij}=a_{ii}+\sum_{j\ne i}a_{ij}
=2\sum_{j\ne i}a_{ij}=0.$
Hence $A\mathbf1_r=\mathbf0_r$.

For $i\ne j$, quadratic reciprocity gives
$a_{ij}+a_{ji} = \left[\frac{p_j}{p_i}\right] + \left[\frac{p_i}{p_j}\right] =s_is_j.$
On the other hand, the diagonal entries of
$A+A^{\mathsf T}$ are zero, while those of
$S+\mathbf s\mathbf s^{\mathsf T}$ are $s_i+s_i^2=0$, since $s_i\in\F_2$.
Hence $A+A^{\mathsf T} =S+\mathbf s\mathbf s^{\mathsf T}.$
Consequently, $Q=A+S =A^{\mathsf T}+\mathbf s\mathbf s^{\mathsf T}.$
This completes the proof.
\end{proof}

The identities established above allow us to
compare the ranks of $\mathcal A_L$ and
$\mathcal A_L+S$. The following lemma shows
that the difference between these ranks is
determined entirely by the residue class of
$L$ modulo $4$.

\begin{lemma}
\label{lem-odd-D-matrix-rank}
Let $L=p_1\cdots p_r$ be a product of distinct odd primes, and let $A=\mathcal A_L$. With $\mathbf s$, $S$, and $Q=A+S$ defined as in Lemma~\ref{lem:matrix-identities}, we have
$$\operatorname{rank}_{\F_2}(A+S) = \operatorname{rank}_{\F_2}A+\eta,~ \text{ where
} ~\eta=\begin{cases}
     0,&L\equiv1\pmod4,\\
     1,&L\equiv3\pmod4.
    \end{cases}$$
\end{lemma}

\begin{proof}
All matrix operations and dimensions in this proof are over
$\F_2$. 
Since $s_i=1$ if and only if $p_i\equiv3\pmod4$, and $L=\displaystyle\prod_{i=1}^rp_i,$
we have $\eta\equiv\displaystyle\sum_{i=1}^r s_i =\mathbf s^{\mathsf T}\mathbf1_r \pmod2.$ Thus, identifying $\eta\in\{0,1\}$ with its image in $\F_2$, we obtain $\mathbf s^{\mathsf T}\mathbf1_r=\eta.$
By Lemma~\ref{lem:matrix-identities},
$$A \mathbf1_r=\mathbf0_r,\quad A+A^{\mathsf T}=S+\mathbf s\mathbf s^{\mathsf T}, \quad Q=A^{\mathsf T}+\mathbf s\mathbf s^{\mathsf T}.$$
Multiplying the second identity on the right by $\mathbf1_r$ and using $A \mathbf1_r=\mathbf 0_r$ and $S \mathbf1_r=\mathbf s$, we obtain
$A^{\mathsf T}\mathbf1_r
=\mathbf s+\mathbf s(\mathbf s^{\mathsf T}\mathbf1_r)
=(1+\eta)\mathbf s.$

\smallskip
\noindent\textit{\textbf{Case $\eta=0$.}}
Here $A^{\mathsf T}\mathbf{1}_r=\mathbf s$, so
$Q=A^{\mathsf T}+\mathbf s\mathbf s^{\mathsf T}
 =A^{\mathsf T}(I_r+\mathbf{1}_r\mathbf s^{\mathsf T}).$
Moreover, since $\mathbf s^{\mathsf T}\mathbf{1}_r=0$ we have $(I_r+\mathbf{1}_r\mathbf s^{\mathsf T})^2
=I_r+\mathbf{1}_r(\mathbf s^{\mathsf T}\mathbf{1}_r)\mathbf s^{\mathsf T}
=I_r.$
Thus $I_r+\mathbf{1}_r\mathbf s^{\mathsf T}$ is invertible.
Right multiplication by this matrix preserves rank, and hence
$\operatorname{rank}Q =\operatorname{rank}A^{\mathsf T}
=\operatorname{rank}A.$

\smallskip
\noindent\textit{\textbf{Case $\eta=1$.}}
Using $A \mathbf{1}_r=\mathbf{0}_r$ and $\mathbf{1}_r^{\mathsf T}\mathbf s=1$, we find $\mathbf{1}_r^{\mathsf T}Q
=\mathbf{1}_r^{\mathsf T}A^{\mathsf T}
 +(\mathbf{1}_r^{\mathsf T}\mathbf s)\mathbf s^{\mathsf T}
=\mathbf s^{\mathsf T}.$
Consequently, if $Q\mathbf x=0$, then $\mathbf s^{\mathsf T}\mathbf x
=\mathbf{1}_r^{\mathsf T}(Q\mathbf x)=0.$
The identity $Q=A^{\mathsf T}+\mathbf s\mathbf s^{\mathsf T}$
then gives $A^{\mathsf T}\mathbf x=\mathbf 0_r$.
Conversely, if $A^{\mathsf T}\mathbf x=\mathbf 0_r$ and
$\mathbf s^{\mathsf T}\mathbf x=0$, then $Q\mathbf x=\mathbf 0_r$.
Therefore $\ker Q=\ker A^{\mathsf T}\cap\mathbf s^\perp,$ where  $\mathbf s^\perp =\{\mathbf x\in\F_2^r:
       \mathbf s^{\mathsf T}\mathbf x=0\}.$
Now $A^{\mathsf T}\mathbf{1}_r=(1+\eta)\mathbf s=\mathbf 0_r$, whereas
$\mathbf s^{\mathsf T}\mathbf{1}_r=1$. Thus the linear functional
\[
\ker A^{\mathsf T}\longrightarrow\F_2,
\qquad
\mathbf x\longmapsto\mathbf s^{\mathsf T}\mathbf x
\]
is nonzero and hence has rank one. Its kernel is exactly
$\ker Q$, so the rank--nullity theorem yields $\dim\ker Q=\dim\ker A^{\mathsf T}-1.$
Applying rank--nullity to the two $r\times r$ matrices gives
$\operatorname{rank}Q
=\operatorname{rank}A^{\mathsf T}+1
=\operatorname{rank}A+1.$
This proves the asserted formula in both cases.
\end{proof}

We next record a general rank identity for
augmented matrices. This identity will be
useful when translating the linear conditions
arising from descent into equivalent matrix
systems and comparing their kernel dimensions.

\begin{lemma}
\label{lem-augmented-matrix-rank}
Let $K$ be a field, let $X\in M_{m,n}(K)$, and let
$\mathbf u\in K^m$ and $\mathbf v\in K^n$ be column vectors.
Then
\[
\operatorname{rank}
[X+\mathbf u\mathbf v^{\mathsf T}\mid\mathbf u]
=
\operatorname{rank}[X\mid\mathbf u]
=
\operatorname{rank}
\begin{pmatrix}
X^{\mathsf T}\\
\mathbf u^{\mathsf T}
\end{pmatrix}.
\]
In particular, if $m=n$ and the common rank is $\rho$, then
$\dim\ker[X+\mathbf u\mathbf v^{\mathsf T}\mid\mathbf u] =n+1-\rho,$
whereas $\dim\ker
\begin{pmatrix}
X^{\mathsf T}\\
\mathbf u^{\mathsf T}
\end{pmatrix}
=n-\rho.$
\end{lemma}

\begin{proof}
Write $\mathbf v=(v_1,\ldots,v_n)^{\mathsf T}$.
For each $1\le j\le n$, subtract $v_j$ times the last column from the $j$th column of $[X+\mathbf u\mathbf v^{\mathsf T}
\mid\mathbf u]$.
These elementary column operations transform this matrix
into $[X\mid\mathbf u]$, so the two matrices have the same
rank.
Moreover, $[X\mid\mathbf u]^{\mathsf T} =
\begin{pmatrix}
X^{\mathsf T}\\
\mathbf u^{\mathsf T}
\end{pmatrix}.$
Since transposition preserves rank, the asserted rank
identities follow.
When $m=n$, the kernel dimensions follow from the
rank--nullity theorem.
\end{proof}

% The preceding lemmas provide the linear-algebraic tools needed to analyze the matrices associated with our descent
% calculations.
We now recall the $2$-isogeny descent method, which relates the Mordell--Weil rank of an elliptic curve
to the images of certain descent homomorphisms. We then describe how these images can be studied through quartic equation.

\subsection{\texorpdfstring{$2$-Descent via $2$-Isogeny}{lg}}

Suppose that $E : Y^2 = X^3 + aX^2 + bX$ is an elliptic curve over $\Q$ where $b\neq 0$ and $a^2-4b\neq 0$. Let $E' : Y^2 = X^3 - 2aX^2 + (a^2 - 4b)X$ be the corresponding $2$-isogenous elliptic curve to $E$ and the corresponding isogeny map for $x\neq 0$ is given by 
\begin{equation}\label{2 isogeny}
    \Phi : E \longrightarrow E', \qquad \Phi(x,y)=\left(\frac{y^2}{x^2},\, \frac{y(x^2-b)}{x^2}\right)
\end{equation}
 which is of degree $2$. At the remaining two points, set
$\Phi(\mathcal{O}_E)=\Phi((0,0))=\mathcal{O}_{E'}$, where $\mathcal{O}_E$ and $\mathcal{O}_{E'}$ are the points at infinity on $E$ and $E'$, respectively.
Let $\Q^{\times }$ denote the multiplicative group of nonzero rational numbers,
and let $\Q^{\times 2}$ be its subgroup consisting of rational squares.
Thus, $\Q^{\times }/\Q^{\times 2}$ denotes the group of nonzero rational numbers
modulo squares. We denote by $E(\Q)$ and
$E'(\Q)$ the groups of
$\Q$-rational points on $E$
and $E'$, respectively. For $d\in\Q^\times$, let $[d]=d(\Q^\times)^2$ denote its square class. We
define the $2$-descent homomorphism $\alpha:E(\Q)\longrightarrow \Q^{\times}/\Q^{\times2}$ by
\begin{equation} \label{alpha}
 \alpha(P)=
\begin{cases}
{[1]}  & \text{if } P=\mathcal{O}_E, \text{ the point at infinity},\\[4pt]
{[b]}  & \text{if } P=(0,0),\\[4pt]
{[x]}  & \text{if } P=(x,y)\text{ with }x\neq 0.
\end{cases}   
\end{equation}
Similarly, we define the $2$-descent homomorphism on the $2$-isogenous
curve $E'$ by $
\alpha':E'(\Q)\longrightarrow
\Q^{\times }/\Q^{\times 2}$
as follows:
\begin{equation}\label{alpha'}
\alpha'(P)=
\begin{cases}
{[1]}  & \text{if } P=\mathcal{O}_{E'}, \text{ the point at infinity},\\[4pt]
{[b']}  & \text{if } P=(0,0),\\[4pt]
{[x]}  & \text{if } P=(x,y)\text{ with }x\neq 0,
\end{cases}
\end{equation}
where $b'= a^2-4b$ is the coefficient of $X$ in the equation of $E'$.

The descent maps \eqref{alpha} and \eqref{alpha'} reduce the problem of determining the rank of an elliptic curve over $\Q$ to the study of their images in $\Q^\times/\Q^{\times 2}$. The following standard results provide a method for computing the rank (see \cite{silverman2009arithmetic} for details).

\begin{proposition}\label{rank proposition}
Let $E/\Q$ be an elliptic curve admitting a $2$-isogeny $\Phi:E\longrightarrow E'$ defined over $\Q$ as in \eqref{2 isogeny}. Let
$r=\operatorname{rank}E(\Q)$. Then 
$$2^{r+2}
=\lvert\operatorname{Im}\alpha\rvert\,\lvert\operatorname{Im}\alpha'\rvert,$$
where $\alpha$ and
$\alpha^\prime$
are the descent maps defined in \eqref{alpha} and
\eqref{alpha'} respectively.
\end{proposition}

By Proposition \ref{rank proposition}, the rank of an elliptic curve
$E$ of the given form can be determined by computing the cardinalities
of the images of the descent maps. The following result provides
a criterion for determining whether an element belongs to these
images. For a proof, see \cite{silverman2009arithmetic}.

\begin{lemma}
\label{image-torsors}
Let $E:Y^2=X^3+aX^2+bX$ and $E':Y^2=X^3-2aX^2+(a^2-4b)X$ be the corresponding $2$-isogenous elliptic curves, where
$a,b\in\mathbb Z$ and $b(a^2-4b)\neq0$.
Use the square-class notation $[u]=u(\Q^\times)^2$
introduced above. Then the following hold:
\begin{itemize}
    \item Every element of $\operatorname{Im}(\alpha)$
     has a representative $d\in\mathbb Z$ such that $d$ is square-free and $d\mid b$. For each such $d$, we have
    \[
    [d]\in\operatorname{Im}(\alpha)
    \quad\Longleftrightarrow\quad
   N^2=du^4+au^2v^2+\frac{b}{d}v^4
    \]
    for some $(N,u,v)\in\mathbb Z^3$ satisfying
    $(u,v)\neq(0,0)$ and $\gcd(u,v)=1$.

    \item Every element of $\operatorname{Im}(\alpha')$
       has a representative $d\in\mathbb Z$ such that $d$ is square-free and
$d\mid a^2-4b$. For each such $d$, we have
    \[
    [d]\in\operatorname{Im}(\alpha')
    \quad\Longleftrightarrow\quad
    N^2=du^4-2au^2v^2+
    \frac{a^2-4b}{d}v^4
    \]
    for some $(N,u,v)\in\mathbb Z^3$ satisfying
    $(u,v)\neq(0,0)$ and $\gcd(u,v)=1$.
\end{itemize}
\end{lemma}

We now specialize to the family of
elliptic curves $E_{\pm L}:y^2=x^3\pm Lx,$ where $L$ is a positive square-free integer. By Lemma~\ref{image-torsors},
determining the descent images amounts to studying the solvability of the corresponding quartic equation. Before doing so, we record some relations between the square classes occurring in these images. These relations allow us to reduce
the number of square-free divisors that must be considered.

\begin{lemma}\label{lem:negative-divisor-equivalence}
Let $L>0$ be a square-free integer,
and let $E_{-L}:y^2=x^3-Lx.$
Let $\alpha$ be the descent map
defined in \eqref{alpha}. Then,
for every positive divisor $\ell$ of
$L$, we have
$$[\ell]\in\operatorname{Im}(\alpha) \iff
[-L/\ell]\in\operatorname{Im}(\alpha), ~\text{ and }~
[-\ell]\in\operatorname{Im}(\alpha) \iff [L/\ell]\in\operatorname{Im}(\alpha).$$
Consequently, the positive and negative
square classes occurring in
$\operatorname{Im}(\alpha)$ determine
one another.
\end{lemma}

\begin{proof}
Since $\alpha$ is a group homomorphism, its image  $H=\operatorname{Im}(\alpha)$ is a subgroup of $\Q^\times/\Q^{\times2}$.
Moreover, $\alpha((0,0))=[-L]\in H.$
Thus, for every positive divisor
$\ell$ of $L$, we have $$[\ell]\in H \iff [-L][\ell]\in H \iff [-L/\ell]\in H,$$
where the last equivalence follows
from $[-L][\ell]=[-L\ell]=[-L/\ell],$ since $\ell^2$ is a rational square.
Similarly,
$$[-\ell]\in H \iff [-L][-\ell]\in H \iff [L/\ell]\in H.$$
This proves the assertions.
\end{proof}

For the curve $E_L$, the positivity of the $x$-coordinates of its nontrivial rational points leads to a simpler description of the descent image. In this case, all square classes in the image are positive, and they occur in complementary pairs, as shown below.

\begin{lemma}\label{lem:positive-divisor-equivalence}
Let $L>0$ be a square-free integer,
and let $E_L:y^2=x^3+Lx.$ Let $\alpha$ be the descent map
defined in \eqref{alpha}. Then every element of $\operatorname{Im}(\alpha)$ has a positive square-free
representative. Moreover, for every positive divisor $\ell$ of $L$, $$[\ell]\in\operatorname{Im}(\alpha) \iff [L/\ell]\in\operatorname{Im}(\alpha).$$
Consequently, it suffices to consider one divisor from each
complementary pair $\{\ell,L/\ell\}$.
\end{lemma}

\begin{proof}
For any rational point $P=(x,y)\in E_L(\Q)$ with
$x\ne0$, we have $y^2=x(x^2+L).$ Since $L>0$, we have $x^2+L>0$, and therefore $x>0$.
Also, $\alpha((0,0))=[L],$ which is represented by a positive
number. Hence every square class in $\operatorname{Im}(\alpha)$ has a positive representative. Since $\alpha$ is a group homomorphism and $[L]\in\operatorname{Im}(\alpha)$,
we obtain 
$$[\ell]\in\operatorname{Im}(\alpha) \iff [L][\ell]\in\operatorname{Im}(\alpha) \iff [L/\ell]\in\operatorname{Im}(\alpha),$$
where the last equivalence follows from $[L\ell]=[L/\ell]$.
\end{proof}

The preceding lemmas reduce the number of square classes that need to be examined in the descent calculations. We now translate necessary local solvability conditions
for the associated quartic into linear equations over
$\F_2$. The following lemma shows that, for a class of quartic equations whose coefficients involve odd powers of the prime divisors of $L$, these conditions can be
expressed using the matrix $\mathcal A_L$.

\begin{lemma}
\label{lem-odd-exponent-descent-conditions}
Let $L=p_1\cdots p_r$ be a product of distinct odd primes,
and let $H=\displaystyle\prod_{i=1}^r p_i^{e_i}$, where $e_i\in\{1,3\}$.
For a positive divisor $\delta$ of $L$, write
$\mathbf x=\varphi(\delta)=(x_1,\ldots,x_r)^{\mathsf T}$ as defined in \eqref{varphi-def}.
Let $\lambda,\mu\in\mathbb Z$ satisfy $\gcd(\lambda\mu,L)=1$.
Suppose that $N^2=\lambda\delta u^4+\mu\frac{H}{\delta}v^4$
has an integer solution $(N,u,v)$ with $\gcd(u,v)=1$, then, for each $1\le i\le r$, we have
\[
(\mathcal A_L\mathbf x)_i=
\begin{cases}
\displaystyle\left[\frac{\lambda}{p_i}\right],&x_i=0,\\[6pt]
\displaystyle\left[\frac{\mu}{p_i}\right],&x_i=1.
\end{cases}
\]
\end{lemma}

\begin{proof}
Since $e_i\in\{1,3\}$ for every $1\le i\le r$,
we may write $H=Lt^2$, where
$t=\displaystyle\prod_{i=1}^r p_i^{(e_i-1)/2}$.
Thus, the given equation becomes
$N^2=\lambda\delta u^4 +\mu\frac{L}{\delta}(tv^2)^2.$
Set $g=\gcd(u^2,tv^2)$. Since $g^2$ divides
both terms on the right-hand side, we have
$g^2\mid N^2$, and hence $g\mid N$.
Putting $N_0=N/g$, $U_0=u^2/g$, and
$V_0=tv^2/g$, we obtain
\begin{equation}\label{eq:primitive-odd-exponent}
N_0^2=\lambda\delta U_0^2
+\mu\frac{L}{\delta}V_0^2,
\qquad \gcd(U_0,V_0)=1.
\end{equation}

Fix $1\le i\le r$ and put $p=p_i$.
We distinguish two cases.

\smallskip
\noindent\textbf{Case 1: $x_i=0$.}
In this case, $p\nmid\delta$ and
$p\mid L/\delta$. Suppose that $p\mid U_0$.
Since $\gcd(U_0,V_0)=1$, we have
$p\nmid V_0$. As $L$ is square-free and
$\gcd(\lambda\mu,L)=1$, the integer
$\mu L/\delta$ is divisible by $p$
but not by $p^2$. Consequently,
$\lambda\delta U_0^2$ is divisible by
$p^2$, whereas $\mu(L/\delta)V_0^2$
is divisible by $p$ but not by $p^2$.
Equation \eqref{eq:primitive-odd-exponent}
would then imply that $p$ divides $N_0^2$
but $p^2$ does not, which is impossible.
Thus $p\nmid U_0$.

Reducing \eqref{eq:primitive-odd-exponent}
modulo $p$, we obtain
$N_0^2\equiv\lambda\delta U_0^2\pmod p$.
Since $p\nmid\lambda\delta U_0$,
the integer $\lambda\delta$ is a
quadratic residue modulo $p$. Hence
$\left[\frac{\lambda\delta}{p}\right]=0$,
which gives
\begin{equation}\label{eq:odd-exp-case-zero}
\left[\frac{\delta}{p}\right]
=\left[\frac{\lambda}{p}\right].
\end{equation}

\smallskip
\noindent\textbf{Case 2: $x_i=1$.}
Here $p\mid\delta$ and $p\nmid L/\delta$.
Suppose that $p\mid V_0$. Then
$p\nmid U_0$, since $\gcd(U_0,V_0)=1$.
As $L$ is square-free and
$\gcd(\lambda\mu,L)=1$, the integer
$\lambda\delta$ is divisible by $p$
but not by $p^2$. Hence
$\lambda\delta U_0^2$ is divisible by
$p$ but not by $p^2$, whereas
$\mu(L/\delta)V_0^2$ is divisible by
$p^2$. Equation
\eqref{eq:primitive-odd-exponent}
would again imply that $p$ divides $N_0^2$
but $p^2$ does not, a contradiction.
Therefore $p\nmid V_0$.

Reducing \eqref{eq:primitive-odd-exponent}
modulo $p$, we get
$N_0^2\equiv\mu(L/\delta)V_0^2\pmod p$.
Since $p\nmid\mu(L/\delta)V_0$,
the integer $\mu L/\delta$ is a
quadratic residue modulo $p$. Thus
$\left[\frac{\mu(L/\delta)}{p}\right]=0$,
and consequently
\begin{equation}\label{eq:odd-exp-case-one}
\left[\frac{L/\delta}{p}\right]
=\left[\frac{\mu}{p}\right].
\end{equation}

Finally, by the definition of $\mathcal A_L$,
we have $(\mathcal A_L\mathbf x)_i
=\displaystyle\sum_{j\ne i}(x_i+x_j)
\left[\frac{p_j}{p_i}\right].$
If $x_i=0$, then
$(\mathcal A_L\mathbf x)_i
=\left[\frac{\delta}{p_i}\right]$,
which equals $\left[\frac{\lambda}{p_i}\right]$
by \eqref{eq:odd-exp-case-zero}.
Similarly, if $x_i=1$, then
$(\mathcal A_L\mathbf x)_i
=\left[\frac{L/\delta}{p_i}\right]$,
which equals $\left[\frac{\mu}{p_i}\right]$
by \eqref{eq:odd-exp-case-one}.
Therefore,
\[
(\mathcal A_L\mathbf x)_i=
\begin{cases}
\left[\frac{\lambda}{p_i}\right],
&x_i=0,\\[4pt]
\left[\frac{\mu}{p_i}\right],
&x_i=1.
\end{cases}
\]
This completes the proof.
\end{proof}

\section{Proof of the main theorems}\label{Proof of the main theorem}    

\begin{proof}[Proof of Theorem~\ref{Thm Positive D}]
Let $E=E_D:~ y^2=x^3+Dx,$ and  $E': y^2=x^3-4Dx,$ be elliptic curves and let $\alpha,\alpha'$ denote their respective descent maps as defined in \eqref{alpha} and \eqref{alpha'}. Clearly from Table~\ref{tab:root-numbers} it follows that the curve $E_D$ has global root number $-1$. Let $A=\mathcal A_D$ be the matrix as defined in \eqref{matrix A_D}. All matrix operations and kernel dimensions below are over
$\F_2$. For $1\leq i \leq n$, write
$s_i=\left[\frac{-1}{\ell_i}\right],~
S=\operatorname{diag}(s_1,\ldots,s_n), ~
Q=A+S,$ and put
$\eta=
\begin{cases}
0,&D\equiv1\pmod4,\\
1,&D\equiv3\pmod4.
\end{cases}~$
By Lemma~\ref{lem-odd-D-matrix-rank}, we have
$\operatorname{rank}_{\F_2}Q
=\operatorname{rank}_{\F_2}A+\eta$.
Since $k=n-\operatorname{rank}_{\F_2}A$, it follows that $\dim\ker A=k$ and $\dim\ker Q=k-\eta$.

We first prove that $|\operatorname{Im}(\alpha)|\le 2^k$.
By Lemma~\ref{image-torsors}, every class in $\operatorname{Im}(\alpha)$ has a squarefree integer representative $d$, with $d\mid D$,
for which the associated quartic equation $N^2=du^4+\frac Dd\,v^4$ has an integer solution satisfying $\gcd(u,v)=1$. In particular, $u,v$ are not both zero. By Lemma~\ref{lem:positive-divisor-equivalence},
every class in $\operatorname{Im}(\alpha)$
has a positive square-free representative.
Consequently, we may take
$d=\delta$, where $\delta$ is a positive
divisor of $D$.
Write $\delta=\displaystyle\prod_{i=1}^n\ell_i^{x_i}$ and
$\mathbf x=\varphi(\delta)=(x_1,\ldots,x_n)^{\mathsf T}\in \F_2^n$
using \eqref{varphi-def} with $L=D$ and the same
ordering of the primes as in \eqref{matrix A_D}. Put $\delta'=D/\delta$.
The corresponding quartic is $N^2=\delta u^4+\delta'v^4$.
Applying Lemma~\ref{lem-odd-exponent-descent-conditions} with $H=D$ and $\lambda=\mu=1$, we obtain $(A\mathbf x)_i=0$ for every $i$. Hence $A\mathbf x=\mathbf0_n$.
Distinct positive square-free divisors of $D$ represent
distinct square classes. Therefore
$[\delta]\mapsto\varphi(\delta)$ defines an injection
from $\operatorname{Im}(\alpha)$ into $\ker A$.
Consequently,
$$|\operatorname{Im}(\alpha)|\le|\ker A|=2^k.$$

We next prove that $|\operatorname{Im}(\alpha')|\le2^{k+1}$.
The image $\operatorname{Im}(\alpha')$ is a subgroup of
$\Q^\times/(\Q^\times)^2$ and contains
$\alpha'((0,0))=[-4D]=[-D]$.
Multiplication by $[-D]$ therefore gives a bijection between
the positive and negative classes in $\operatorname{Im}(\alpha')$.
By Lemma~\ref{image-torsors}, every positive class in
$\operatorname{Im}(\alpha^\prime)$ has a positive squarefree representative dividing $4D$. Such a representative is either $\delta$ or $2\delta$, where $\delta\mid D$.
Let $J_0$ and $J_1$ denote the sets of positive classes
represented by $\delta$ and $2\delta$, respectively.
Then $|\operatorname{Im}(\alpha^\prime)|=2\bigl(|J_0|+|J_1|\bigr).$

For a class $[\delta]\in J_0$, put
$\mathbf x=\varphi(\delta)$ as in \eqref{varphi-def} and $\delta'=D/\delta$ as above.
The corresponding quartic is
$N^2=\delta u^4-4\delta'v^4$, with $\gcd(u,v)=1$.
Applying Lemma~\ref{lem-odd-exponent-descent-conditions} with $H=D$,
$\lambda=1$, and $\mu=-4$, we find
\[
(A\mathbf x)_i=
\begin{cases}
\left[\frac{-4}{\ell_i}\right]=s_i,&x_i=1,\\
\left[\frac{1}{\ell_i}\right]=0,&x_i=0.
\end{cases}
\]
Thus $A\mathbf x=S\mathbf x$, or equivalently
$Q\mathbf x=\mathbf0_n$.
It follows that $|J_0|\le|\ker Q|=2^{k-\eta}$.

For a class $[2\delta]\in J_1$, put
$\mathbf x=\varphi(\delta)$ as defined in \eqref{varphi-def} and $\delta'=D/\delta$ as above.
The corresponding quartic is
$N^2=2\delta u^4-2\delta'v^4$, with $\gcd(u,v)=1$.
Write $b_i=\bigl[\frac{2}{\ell_i}\bigr]$ and
$\mathbf b=(b_1,\ldots,b_n)^{\mathsf T}$.
Applying Lemma~\ref{lem-odd-exponent-descent-conditions} with $H=D$,
$\lambda=2$, and $\mu=-2$, we obtain
$$(A\mathbf x)_i=
\begin{cases}
\left[\frac{-2}{\ell_i}\right]=b_i+s_i,&x_i=1,\\
\left[\frac{2}{\ell_i}\right]=b_i,&x_i=0.
\end{cases}$$
Hence $Q\mathbf x=\mathbf b$.
If this system has a solution $\mathbf x_0$, its solution set
is $\mathbf x_0+\ker Q$. Otherwise, its solution set is empty.
Consequently, $|J_1|\le2^{k-\eta}$.

Suppose first that $D\equiv3$ or $15\pmod{16}$.
Then $\eta=1$, and the preceding bounds give
$$|\operatorname{Im}(\alpha^\prime)|=2\bigl(|J_0|+|J_1|\bigr)
\le2\bigl(2^{k-1}+2^{k-1}\bigr) =2^{k+1}.$$

Now suppose that $D\equiv5$ or $13\pmod{16}$.
Here $\eta=0$ and $D\equiv5\pmod8$.
We claim that $J_1$ is empty.
Indeed, suppose that $N^2=2\delta u^4-2\delta'v^4$ and $\delta\delta'=D$, with $\gcd(u,v)=1.$ The right-hand side is even, so $N$ is even. Reduction modulo $4$ shows that $u$ and $v$ have the same parity, and primitivity therefore forces both to be odd. Since $\delta$ is odd, $\delta^2\equiv1\pmod8$, and hence
$\delta'\equiv D\delta\equiv 5 \delta\pmod 8.$
Thus $\delta-\delta'\equiv-4\delta\equiv4\pmod8.$
Odd fourth powers are congruent to $1$ modulo $16$.
Reducing the quartic modulo $16$ therefore gives $N^2\equiv2(\delta-\delta')\equiv8\pmod{16},$ which is impossible for a square. Hence $J_1=\varnothing$, and $|\operatorname{Im}(\alpha^\prime)|=2|J_0|\le2^{k+1}.$

Writing $r=\operatorname{rank}E_D(\Q)$,
Proposition~\ref{rank proposition} now gives
$$2^r=\frac{|\operatorname{Im}(\alpha)|\,|\operatorname{Im}(\alpha^\prime)|}{4} \le\frac{2^k\,2^{k+1}}4 =2^{2k-1}.$$
Consequently, $\operatorname{rank}E_D(\Q)\le2k-1.$

Finally, if $\operatorname{rank}_{\F_2}A=n-1$,
then $k=1$, and the rank bound gives
$\operatorname{rank}E_D(\Q)\le1$.
If the Mordell--Weil parity conjecture holds for $E_D$, then
$(-1)^{\operatorname{rank}E_D(\Q)}=\omega(E_D)=-1$.
The rank is therefore odd, and the preceding upper bound
forces $\operatorname{rank}E_D(\Q)=1$.

\end{proof}

\begin{proof}[Proof of Theorem~\ref{Thm Negative D}]

Let $E=E_{-D}:y^2=x^3-Dx$ and $E':y^2=x^3+4Dx$, be elliptic curves and let $\alpha,\alpha'$ be the corresponding descent maps defined as in \eqref{alpha} and \eqref{alpha'} respectively.
The root-number formula in Table~\ref{tab:root-numbers} gives $\omega(E_{-D})=-1$
under the congruence hypotheses of the theorem.
Let $A=\mathcal A_D$ be the matrix defined in
\eqref{matrix A_D}. All matrix operations and kernel
dimensions below are over $\F_2$.
Write $s_i=\bigl[\frac{-1}{\ell_i}\bigr]$,
$S=\operatorname{diag}(s_1,\ldots,s_n)$, and $Q=A+S$,
and put
$\eta=\begin{cases}
0,&D\equiv1\pmod4,\\
1,&D\equiv3\pmod4.
\end{cases}$
By Lemma~\ref{lem-odd-D-matrix-rank}, we have
$\operatorname{rank}_{\F_2}Q
=\operatorname{rank}_{\F_2}A+\eta$.
Since $k=n-\operatorname{rank}_{\F_2}A$, it follows that $\dim\ker A=k$ and $\dim\ker Q=k-\eta$.

We first prove that $|\operatorname{Im}(\alpha)|
\le2^{k+1-\eta}$. Since $\alpha((0,0))=[-D]$, the image
$\operatorname{Im}(\alpha)$ contains $[-D]$.
By Lemma~\ref{lem:negative-divisor-equivalence},
multiplication by $[-D]$ gives a bijection
between the positive and negative square
classes in $\operatorname{Im}(\alpha)$.
Thus it suffices to bound the number
of positive classes. By Lemma~\ref{image-torsors}, every positive class in $\operatorname{Im}(\alpha)$ has a positive square-free representative $\delta\mid D$ for which $N^2=\delta u^4-\delta'v^4$, where $\delta'=D/\delta$,
has an integer solution satisfying $\gcd(u,v)=1$. Write $\delta=\displaystyle\prod_{i=1}^n\ell_i^{x_i}$ and
$\mathbf x=\varphi(\delta)=(x_1,\ldots,x_n)^{\mathsf T}\in \F_2^n$
using \eqref{varphi-def} with $L=D$ and the same
ordering of the primes as in \eqref{matrix A_D}. Applying Lemma~\ref{lem-odd-exponent-descent-conditions} with $H=D$,
$\lambda=1$, and $\mu=-1$, we obtain
\[
(A\mathbf x)_i=
\begin{cases}
\left[\frac{-1}{\ell_i}\right]=s_i,&x_i=1,\\
\left[\frac{1}{\ell_i}\right]=0,&x_i=0.
\end{cases}
\]
Thus $A\mathbf x=S\mathbf x$, or equivalently
$Q\mathbf x=\mathbf0_n$. Distinct positive square-free divisors of $D$ represent
distinct square classes. Hence
$[\delta]\mapsto\varphi(\delta)$ injects the positive
classes in $\operatorname{Im}(\alpha)$ into $\ker Q$.
Therefore
$|\operatorname{Im}(\alpha)|\le2|\ker Q|=2^{k+1-\eta}$.

We next bound $|\operatorname{Im}(\alpha')|$.
By Lemma~\ref{image-torsors}, every class in this image has
a signed square-free integer representative $d\mid4D$
for which $N^2=du^4+\frac{4D}{d}v^4$ has an integer
solution satisfying $\gcd(u,v)=1$.
If $d<0$, then both $d$ and $4D/d$ are negative, so the
right-hand side is strictly negative, a contradiction.
Thus every class in $\operatorname{Im}(\alpha')$ is positive
and has a representative either $\delta$ or $2\delta$,
where $\delta>0$ and $\delta\mid D$.

Let $J_0$ and $J_1$ denote the sets of classes represented
by $\delta$ and $2\delta$, respectively.
These sets are disjoint, and $|\operatorname{Im}(\alpha')|=|J_0|+|J_1|.$

For a class $[\delta]\in J_0$, put
$\mathbf x=\varphi(\delta)$ as in \eqref{varphi-def} and $\delta'=D/\delta$ as above.
The corresponding quartic is
$N^2=\delta u^4+4\delta'v^4$, with $\gcd(u,v)=1$.
Applying Lemma~\ref{lem-odd-exponent-descent-conditions} with $H=D$,
$\lambda=1$, and $\mu=4$, we obtain
\[
(A\mathbf x)_i=
\begin{cases}
\left[\frac{4}{\ell_i}\right]=0,&x_i=1,\\
\left[\frac{1}{\ell_i}\right]=0,&x_i=0.
\end{cases}
\]
Hence $A\mathbf x=\mathbf0_n$, and consequently
$|J_0|\le|\ker A|=2^k$.

For a class $[2\delta]\in J_1$, put
$\mathbf x=\varphi(\delta)$ as defined in \eqref{varphi-def} and $\delta'=D/\delta$ as above.
The corresponding quartic is
$N^2=2\delta u^4+2\delta'v^4$, with $\gcd(u,v)=1$.
Write $b_i=\bigl[\frac{2}{\ell_i}\bigr]$ and
$\mathbf b=(b_1,\ldots,b_n)^{\mathsf T}$.
Applying Lemma~\ref{lem-odd-exponent-descent-conditions} with $H=D$
and $\lambda=\mu=2$, we obtain
$(A\mathbf x)_i=b_i$ for every $i$.
Thus $A\mathbf x=\mathbf b$.
If this system has a solution $\mathbf x_0$, its solution
set is $\mathbf x_0+\ker A$. Otherwise, its solution set
is empty. Consequently, $|J_1|\le2^k$, and hence
$|\operatorname{Im}(\alpha')|\le2^{k+1}$.

Suppose first that $D\equiv7$ or $15\pmod{16}$.
Then in this case we have $\eta=1$, so the preceding bounds give $|\operatorname{Im}(\alpha)|\le2^k,~ |\operatorname{Im}(\alpha')|\le2^{k+1}.$
Now suppose that $D\equiv5$ or $9\pmod{16}$. Here $\eta=0$.
We claim that $J_1$ is empty. Indeed, suppose that
$N^2=2\delta u^4+2\delta'v^4, ~ \delta\delta'=D, ~ \gcd(u,v)=1. $
The right-hand side is even, so $N$ is even.
Reduction modulo $4$ shows that $u$ and $v$ have the
same parity, and primitivity therefore forces both to
be odd.
Since odd fourth powers are congruent to $1$ modulo $16$,
reduction modulo $32$ gives $N^2\equiv2(\delta+\delta')\pmod{32}.$
An even square is congruent to $0$, $4$, or $16$ modulo $32$.
Therefore $\delta+\delta'\equiv0,2,\text{ or }8\pmod{16}.$
For odd $\delta$ and $\delta'$, a direct calculation
modulo $16$ shows that $\delta+\delta'\equiv0\text{ or }8\pmod{16}$ which implies $\delta\delta'\equiv7\text{ or }15\pmod{16},$ whereas $\delta+\delta'\equiv2\pmod{16}$ implies $\delta\delta'\equiv1\text{ or }13\pmod{16}.$
Each possibility contradicts $D\equiv 5$ or $9\pmod{16}$.
Thus $J_1=\varnothing$, and hence
$|\operatorname{Im}(\alpha)|\le2^{k+1}$ and $|\operatorname{Im}(\alpha')|=|J_0|\le2^k.$
Writing $r=\operatorname{rank}E_{-D}(\Q)$,
Proposition~\ref{rank proposition} now gives, in either case,
$$2^r =\frac{|\operatorname{Im}(\alpha)|\,
        |\operatorname{Im}(\alpha')|}{4}
\le\frac{2^{2k+1}}4 =2^{2k-1}.$$
Consequently, we obtain $\operatorname{rank}E_{-D}(\Q)\le2k-1.$

Finally, if $\operatorname{rank}_{\F_2}A=n-1$,
then $k=1$, and the rank bound gives
$\operatorname{rank}E_{-D}(\Q)\le1$.
If the Mordell--Weil parity conjecture holds for $E_{-D}$,
then
$(-1)^{\operatorname{rank}E_{-D}(\Q)}
=\omega(E_{-D})=-1$.
The rank is therefore odd, and the preceding upper bound
forces $\operatorname{rank}E_{-D}(\Q)=1$.
\end{proof}

\begin{proof}[Proof of Theorem~\ref{Thm negative 2D}]
Let $E=E_{-2D}:y^2=x^3-2Dx$ and $E':y^2=x^3+8Dx$,
and let $\alpha,\alpha'$ be the corresponding descent maps as defined in \eqref{alpha} and \eqref{alpha'} respectively.
The root-number formula in Table~\ref{tab:root-numbers} gives $\omega(E)=-1$.
Let $A=\mathcal A_D$ be the matrix defined in
\eqref{matrix A_D}.
All matrix operations and kernel dimensions below are
over $\F_2$.
For $1\le i\le n$, put
$b_i=\bigl[\frac{2}{\ell_i}\bigr]$ and
$s_i=\bigl[\frac{-1}{\ell_i}\bigr]$.
Write $\mathbf s=(s_1,\ldots,s_n)^{\mathsf T}$,
$B=\operatorname{diag}(b_1,\ldots,b_n)$, and
$S=\operatorname{diag}(s_1,\ldots,s_n)$.
Put $C=A+B$, $T=A+B+S$, and $\kappa=n-\rho_D$.
Thus $\mathcal R_D=[T\mid\mathbf s]$, as in
\eqref{eq:augmented-RD}, and
$\dim\ker\mathcal R_D=\kappa+1$.

We first establish a matrix-rank identity that will be used
to bound the descent image.
By Lemma~\ref{lem:matrix-identities},
$A+A^{\mathsf T}=S+\mathbf s\mathbf s^{\mathsf T}$.
Since $B$ is symmetric, it follows that
$T=C^{\mathsf T}+\mathbf s\mathbf s^{\mathsf T}$.
Applying Lemma~\ref{lem-augmented-matrix-rank} with
$X=C^{\mathsf T}$ and $\mathbf u=\mathbf v=\mathbf s$,
we obtain
\begin{equation}
\label{eq-negative-2D-rank-identity}
\operatorname{rank}_{\F_2}
\begin{pmatrix}
C\\ \mathbf s^{\mathsf T}
\end{pmatrix}
=\operatorname{rank}_{\F_2}\mathcal R_D =\rho_D.
\end{equation}
Consequently, the rank--nullity theorem gives
$\dim\ker\mathcal R_D=\kappa+1$ and
$\dim\ker\binom{C}{\mathbf s^{\mathsf T}}=\kappa$.

We next establish the bound $|\operatorname{Im}(\alpha)|
\le2^{\kappa+2}$.
By Lemma~\ref{image-torsors}, every class in
$\operatorname{Im}(\alpha)$ has a signed square-free integer
representative supported on the primes dividing $2D$.
Moreover, $\operatorname{Im}(\alpha)$ is a subgroup of
$\Q^\times/(\Q^\times)^2$ and contains
$\alpha((0,0))=[-2D]$. Hence multiplication by
$[-2D]$ gives a bijection between the classes with odd
and even square-free representatives.

Let $I_0$ denote the set of classes in
$\operatorname{Im}(\alpha)$ with odd squarefree representatives. Then $|\operatorname{Im}(\alpha)|=2|I_0|.$
Every class in $I_0$ has a unique representative
$d=(-1)^t\delta$, where $t\in\{0,1\}$ and
$\delta>0$ with $\delta\mid D$.
Write $\delta=\displaystyle\prod_{i=1}^n\ell_i^{x_i}$ and
$\mathbf x=\varphi(\delta)=(x_1,\ldots,x_n)^{\mathsf T}$,
using \eqref{varphi-def} with $L=D$ and the same
prime order as in \eqref{matrix A_D}.
Put $\delta'=D/\delta$.

By Lemma~\ref{image-torsors}, the corresponding quartic
$N^2=(-1)^t\delta u^4+(-1)^{t+1}2\delta'v^4$
has an integer solution satisfying $\gcd(u,v)=1$.
Applying Lemma~\ref{lem-odd-exponent-descent-conditions} with $H=D$,
$\lambda=(-1)^t$, and $\mu=(-1)^{t+1}2$, we obtain
\[
(A\mathbf x)_i=
\begin{cases}
\left[\frac{(-1)^{t+1}2}{\ell_i}\right]
=b_i+(t+1)s_i,&x_i=1,\\
\left[\frac{(-1)^t}{\ell_i}\right]
=ts_i,&x_i=0.
\end{cases}
\]
Thus $A\mathbf x=(B+S)\mathbf x+t\mathbf s$,
or equivalently $T\mathbf x+t\mathbf s=\mathbf0_n$. Distinct signed square-free divisors of $D$ represent
distinct square classes. Hence the map
$[(-1)^t\delta]\mapsto\binom{\varphi(\delta)}{t}$,
where $t \in \F_2$,
embeds $I_0$ into $\ker\mathcal R_D$.
Consequently, $|I_0|\le2^{\kappa+1}$, and therefore
$|\operatorname{Im}(\alpha)|=2|I_0|\le2^{\kappa+2}$.

We next prove that $|\operatorname{Im}(\alpha')|
\le2^{\kappa+1}$.
By Lemma~\ref{image-torsors}, every class in
$\operatorname{Im}(\alpha')$ has a signed square-free integer
representative $d\mid8D$ for which
$N^2=du^4+\frac{8D}{d}v^4$ has an integer solution
satisfying $\gcd(u,v)=1$.
If $d<0$, then both $d$ and $8D/d$ are negative, so the
right-hand side is strictly negative, a contradiction.
Thus every class in $\operatorname{Im}(\alpha')$ is positive
and has a representative either $\delta$ or $2\delta$,
where $\delta>0$ and $\delta\mid D$.

Let $J_0$ and $J_1$ denote the sets of classes represented
by $\delta$ and $2\delta$, respectively.
The image $\operatorname{Im}(\alpha')$ contains
$\alpha'((0,0))=[8D]=[2D]$.
Multiplication by $[2D]$ gives a bijection between
$J_0$ and $J_1$. Therefore, 
$|\operatorname{Im}(\alpha')|=2|J_0|$.

For a class $[\delta]\in J_0$, put
$\mathbf x=\varphi(\delta)$ as defined in \eqref{varphi-def} and $\delta'=D/\delta$ as above.
The corresponding quartic is
$N^2=\delta u^4+8\delta'v^4$, with $\gcd(u,v)=1$.
Applying Lemma~\ref{lem-odd-exponent-descent-conditions} with $H=D$, $\lambda=1$, and $\mu=8$, we obtain
\[
(A\mathbf x)_i=
\begin{cases}
\left[\frac{8}{\ell_i}\right]=b_i,&x_i=1,\\
\left[\frac{1}{\ell_i}\right]=0,&x_i=0.
\end{cases}
\]
Hence $A\mathbf x=B\mathbf x$, or equivalently
$C\mathbf x=\mathbf0_n$. If $u$ were even, the condition $\gcd(u,v)=1$ would force $v$ to be odd. Since $\delta$ and
$D/\delta$ are odd, reduction modulo $16$ would give
$N^2=\delta u^4+8\frac{D}{\delta}v^4 \equiv8\pmod{16}.$
This is impossible, since a square is congruent to
$0$, $1$, $4$, or $9$ modulo $16$.
Hence $u$ is odd, and the quartic equation shows
that $N$ is odd. Reducing the equation modulo $8$
now gives $\delta\equiv1\pmod8$, and consequently
$\delta\equiv1\pmod4$.
Since $s_i=1$ precisely when $\ell_i\equiv3\pmod4$,
the factorisation $\delta=\prod_{i=1}^n\ell_i^{x_i}$
shows that $\sum_{i=1}^n s_ix_i$ is even.
Consequently,
$\mathbf s^{\mathsf T}\mathbf x
=\sum_{i=1}^n s_ix_i=0$ in $\F_2$.
It follows that the map $[\delta]\mapsto\varphi(\delta)$ embeds
$J_0$ into $\ker\binom{C}{\mathbf s^{\mathsf T}}$.
By \eqref{eq-negative-2D-rank-identity} and the
rank--nullity theorem, this kernel has dimension $\kappa$.
Hence $|J_0|\le2^\kappa$, and therefore
$|\operatorname{Im}(\alpha')|=2|J_0|\le2^{\kappa+1}$.
Writing $r=\operatorname{rank}E_{-2D}(\Q)$,
Proposition~\ref{rank proposition} gives
\begin{equation}
\label{eq-negative-2D-descent-rank}
2^r=\frac{|\operatorname{Im}(\alpha)|\,
      |\operatorname{Im}(\alpha')|}{4}.
\end{equation}
The preceding bounds give
$2^r\le2^{\kappa+2}2^{\kappa+1}/4=2^{2\kappa+1}.$
Consequently, $r\le2\kappa+1=2(n-\rho_D)+1$.

If the Mordell--Weil parity conjecture holds for $E$,
then $(-1)^r=\omega(E)=-1$.
Thus $r$ is odd, and the upper bound gives
$r\in\{1,3,\ldots,2(n-\rho_D)+1\}$.
In particular, if $\rho_D=n$, then $\kappa=0$ and
$r\le1$ unconditionally. Under the parity conjecture,
$r$ is odd and therefore equals one.
This proves the second assertion.

Finally, suppose that the hypotheses of the third assertion
hold, and put $d=\varepsilon\delta.$
Since $\varepsilon^{-1}=\varepsilon$, the given square
condition becomes $m^2=d-\frac{2D}{d}h^4.$
Multiplying by $d^2/h^6$, we obtain
$$\left(\frac{dm}{h^3}\right)^2 =
\left(\frac d{h^2}\right)^3
-2D\left(\frac d{h^2}\right).$$
Hence $P=(d/h^2,dm/h^3)\in E(\Q)$, and
$\alpha(P)=[d]$ by the definition of $\alpha$.
The integer $d$ is odd and square-free. If $d=1$, the square condition gives $m^2=1-2Dh^4<0$, a contradiction. Hence $d\ne1$, and therefore $[d]\ne[1]$.
Moreover, $[d]\ne[-2D]$. Indeed, equality would imply
that $-2Dd$ is a rational square and hence an integer
square. But $d$ and $D$ are odd, so
$-2Dd\equiv2\pmod4$, which is impossible for a square.
Thus $[d]$ and $[-2D]$ generate a subgroup of order four
in $\operatorname{Im}(\alpha)$, namely
$\{[1],[d],[-2D],[-2Dd]\}$.
Therefore $|\operatorname{Im}(\alpha)|\ge4$.

Also, $[1]$ and $[2D]$ are distinct classes in
$\operatorname{Im}(\alpha')$, so
$|\operatorname{Im}(\alpha')|\ge2$.
Equation~\eqref{eq-negative-2D-descent-rank} now gives
$$2^r=\frac{|\operatorname{Im}(\alpha)|\,
      |\operatorname{Im}(\alpha')|}{4}
\ge\frac{4\cdot2}{4}=2.$$
Thus $r\ge1$.
Since $\rho_D=n$ already gives $r\le1$, we conclude that $\operatorname{rank}E_{-2D}(\Q)=1$
without assuming the parity conjecture.
\end{proof}

\begin{corollary}[The two-prime specialization]
\label{cor-two-prime-specialization}
Let $p$ and $q$ be distinct odd primes, and put
$E=E_{-2pq}:y^2=x^3-2pqx$,
$r=\operatorname{rank}E(\Q)$, and
$\chi=\left(\frac{p}{q}\right)$.
Then $\omega(E)=-1$. Moreover, assuming the Mordell--Weil parity conjecture for $E$, the following conclusions hold.
\begin{enumerate}
\item If $(p,q)\equiv(1,1)\pmod8$, then
$r\in\{1,3,5\}$ when $\chi=1$, and
$r\in\{1,3\}$ when $\chi=-1$.

\item If one prime is congruent to $1\pmod8$ and the other
is congruent to $3$, $5$, or $7\pmod8$, then
$r\in\{1,3\}$ when $\chi=1$, and $r=1$ when $\chi=-1$.

\item If, up to interchange,
$(p,q)\equiv(3,3),(3,5),(3,7),(5,5)$, or $(5,7)\pmod8$,
then $r=1$ for either value of $\chi$.

\item If $(p,q)\equiv(7,7)\pmod8$, then
$r\in\{1,3\}$ for either value of $\chi$.
\end{enumerate}
The upper bounds for $r$ in these conclusions hold
unconditionally. Moreover, the square conditions in
\cite[Theorem~1.5]{Arkabrata-Bidisha-Richa}
give $r=1$ without assuming the parity conjecture.
\end{corollary}

\begin{proof}
Apply Theorem~\ref{Thm negative 2D} with $D=pq$ and $n=2$.
It gives $\omega(E)=-1$ and $r\le5-2\rho_{pq}$, where
$\rho_{pq}=\operatorname{rank}_{\F_2}\mathcal R_{pq}$.
We compute $\rho_{pq}$ explicitly.
Put $t=\bigl[\frac{p}{q}\bigr]$,
$s_p=\bigl[\frac{-1}{p}\bigr]$,
$s_q=\bigl[\frac{-1}{q}\bigr]$,
$b_p=\bigl[\frac{2}{p}\bigr]$, and
$b_q=\bigl[\frac{2}{q}\bigr]$.
All matrix entries and operations below are over
$\F_2$.
By the definition of the additive Legendre symbol,
$\chi=(-1)^t$.
Quadratic reciprocity gives
$\bigl[\frac{q}{p}\bigr]=t+s_ps_q$.
Therefore, using \eqref{matrix A_D} and
\eqref{eq:augmented-RD} in the prime order $(p,q)$, we obtain
\[
\mathcal A_{pq}
=
\begin{pmatrix}
t+s_ps_q&t+s_ps_q\\
t&t
\end{pmatrix},
\qquad
\mathcal R_{pq}
=
\begin{pmatrix}
t+s_ps_q+b_p+s_p&t+s_ps_q&s_p\\
t&t+b_q+s_q&s_q
\end{pmatrix}.
\]

For an odd prime $\ell$, the pair
$\bigl([\frac{-1}{\ell}],[\frac{2}{\ell}]\bigr)$ equals
$(0,0)$, $(1,1)$, $(0,1)$, or $(1,0)$ according as
$\ell\equiv1,3,5$, or $7\pmod8$.
Substitution into $\mathcal R_{pq}$ and row reduction give
\[
\begin{array}{c|cc}
(p,q)\bmod8
&\rho_{pq}\text{ when }\chi=1
&\rho_{pq}\text{ when }\chi=-1\\ \hline
(1,1)&0&1\\
(1,3),(1,5),(1,7)&1&2\\
(3,3),(3,5),(3,7),(5,5),(5,7)&2&2\\
(7,7)&1&1
\end{array}
\]
Here the residue pairs are listed up to interchange.

For example, if $(p,q)\equiv(1,1)\pmod8$, then
$\mathcal R_{pq}=
\begin{pmatrix}
t&t&0\\
t&t&0
\end{pmatrix}$,
which has rank zero when $t=0$ and rank one when $t=1$.
If $(p,q)\equiv(7,7)\pmod8$, then
$\mathcal R_{pq}=
\begin{pmatrix}
t&t+1&1\\
t&t+1&1
\end{pmatrix}$,
which has rank one for either value of $t$.
The remaining substitutions give the middle two rows
of the table.

Interchanging $p$ and $q$ permutes the rows and the first
two columns of $\mathcal R_{pq}$, preserving its rank.
Whenever the tabulated rank depends on $\chi$, at least
one prime is congruent to $1\pmod8$.
Quadratic reciprocity then gives
$\left(\frac{p}{q}\right)=\left(\frac{q}{p}\right)$.
Thus the table is well defined for residue pairs listed
up to interchange.

The unconditional bound $r\le5-2\rho_{pq}$ gives
$r\le5$, $r\le3$, or $r\le1$ according as
$\rho_{pq}=0$, $1$, or $2$.
Under the Mordell--Weil parity conjecture,
$(-1)^r=\omega(E)=-1$, so $r$ is odd.
Consequently,
\[
r\in
\begin{cases}
\{1,3,5\},&\rho_{pq}=0,\\
\{1,3\},&\rho_{pq}=1,\\
\{1\},&\rho_{pq}=2.
\end{cases}
\]
Combining these possibilities with the rank table proves
the four assertions and their unconditional upper bounds.
This recovers
\cite[Theorems~1.1--1.4]{Arkabrata-Bidisha-Richa}.
\end{proof}

\begin{proof}[Recovery of the explicit rank-one criteria]
The preceding rank calculation gives $\rho_{pq}=2$
for every residue pair in the table below.
Hence Theorem~\ref{Thm negative 2D} gives $r\le1$
unconditionally.

Part~\textup{(3)} of Theorem~\ref{Thm negative 2D} gives $r=1$ whenever
there exist a positive divisor $\delta\mid pq$,
a sign $\varepsilon\in\{1,-1\}$, and positive integers
$h,m$ such that $m^2=\varepsilon\left(\delta-2h^4\frac{pq}{\delta}\right)$.
The square conditions in
\cite[Theorem~1.5]{Arkabrata-Bidisha-Richa}
are obtained from the following substitutions.
The residue pairs in this table are ordered, and
$\chi$ always denotes $\left(\frac{p}{q}\right)$.

\begingroup
\setlength{\arraycolsep}{5pt}
\small
\[
\begin{array}{c|c|c|ccc}
(p,q)\bmod 8 & \chi & m^2 & \delta & \varepsilon & h\\ \hline
(3,3)&\pm1&pq-32&pq&1&2\\ \hline
(3,5)&1&p-2q&p&1&1\\
(3,5)&-1&2p-q&q&-1&1\\ \hline
(5,3)&1&q-2p&q&1&1\\
(5,3)&-1&2q-p&p&-1&1\\ \hline
(3,7)&1&32p-q&q&-1&2\\
(3,7)&-1&q-2p&q&1&1 \\ \hline
\end{array}
\qquad
\begin{array}{c|c|c|ccc}
(p,q)\bmod 8 & \chi & m^2 & \delta & \varepsilon & h\\ \hline
(5,5)&\pm1&2pq-1&1&-1&1 \\ \hline
(7,3)&1&p-2q&p&1&1\\
(7,3)&-1&32q-p&p&-1&2\\ \hline
(5,7)&1&32p-q&q&-1&2\\
(5,7)&-1&2q-p&p&-1&1\\ \hline
(7,5)&1&2p-q&q&-1&1\\
(7,5)&-1&32q-p&p&-1&2\\ \hline
\end{array}
\]
\endgroup

In every row, $\delta$ is a positive divisor of $pq$, the expression in the third
column equals
$\varepsilon\left(\delta-2h^4\frac{pq}{\delta}\right)$.
Each expression is odd. Thus, whenever it is a perfect
square, it has a positive integer square root $m$.
Theorem~\ref{Thm negative 2D}\textup{(3)} therefore gives
$\operatorname{rank}E_{-2pq}(\Q)=1$
in every listed case, without assuming the parity conjecture.
This recovers
\cite[Theorem~1.5]{Arkabrata-Bidisha-Richa}.
\end{proof}

\begin{proof}[Proof of Theorem~\ref{thm odd prime exponent extension}]
Throughout, let $H=\prod_{i=1}^n\ell_i^{e_i}$ with
$e_i\in\{1,3\}$, and put $D=\prod_{i=1}^n\ell_i$,
where $n\ge2$ and the $\ell_i$ are distinct odd primes.
Write $H=DF^2$, where $F=\prod_{e_i=3}\ell_i$.
Since $F$ is odd, $H\equiv D\pmod8$.

Let $A=\mathcal A_D$ be the matrix defined in
\eqref{matrix A_D}, and put
$k=n-\operatorname{rank}_{\F_2}A$.
For $1\le i\le n$, write
$b_i=\bigl[\frac{2}{\ell_i}\bigr]$ and
$s_i=\bigl[\frac{-1}{\ell_i}\bigr]$.
Set $\mathbf b=(b_1,\ldots,b_n)^{\mathsf T}$,
$\mathbf s=(s_1,\ldots,s_n)^{\mathsf T}$,
$B=\operatorname{diag}(b_1,\ldots,b_n)$, and
$S=\operatorname{diag}(s_1,\ldots,s_n)$.
Let $\mathcal R_D$ be as in \eqref{eq:augmented-RD},
and put
$\rho_D=\operatorname{rank}_{\F_2}\mathcal R_D$
and $\kappa=n-\rho_D$.

The extended root-number formula in Section~1 applies
to $H$, $-H$, and $-2H$.
In part~(1), $H>0$ and $c_2(H)=-1$.
In part~(2), $-H\equiv11,9,7,1\pmod{16}$ under the
stated hypotheses, so $c_2(-H)=1$.
In part~(3), $-2H<0$ and $v_2(-2H)=1$.
Thus the global root number is $-1$ in each case.

\begin{enumerate}[leftmargin=*]
\item Let $E=E_H: y^2=x^3+Hx,~~ E': y^2=x^3-4Hx,$
and let $\alpha,~\alpha'$ be the corresponding descent maps.
Put $Q=A+S, ~
\eta=\begin{cases}
0,&H\equiv1\pmod4,\\
1,&H\equiv3\pmod4.
\end{cases}$
Since $H\equiv D\pmod4$, Lemma~\ref{lem-odd-D-matrix-rank}
gives $\dim\ker A=k, ~\dim\ker Q=k-\eta.$

We first prove that $|\operatorname{Im}(\alpha)|\le2^k.$
By Lemma~\ref{image-torsors}, every class in
$\operatorname{Im}(\alpha)$ has a signed squarefree
representative supported on the primes dividing $D$.
A negative representative is impossible, since both
coefficients in its associated quartic would be negative.
Thus every class has a positive representative $\delta\mid D$
for which
$N^2=\delta u^4+\frac H\delta v^4, ~ \gcd(u,v)=1.$
Write $\delta=\prod_{i=1}^n\ell_i^{x_i}$ and
$\mathbf x=\varphi(\delta)=(x_1,\ldots,x_n)^{\mathsf T}$,
using \eqref{varphi-def} with $L=D$ and the same
prime order as in \eqref{matrix A_D}.
Applying Lemma~\ref{lem-odd-exponent-descent-conditions}
with $\lambda=\mu=1$, we obtain
$A\mathbf x=\mathbf0_n$.
Distinct positive square-free divisors of $D$ represent
distinct square classes. Hence
$[\delta]\mapsto\varphi(\delta)$ injects
$\operatorname{Im}(\alpha)$ into $\ker A$, and therefore
$|\operatorname{Im}(\alpha)|\le|\ker A|=2^k$.

We next prove that $|\operatorname{Im}(\alpha')|\le2^{k+1}$.
The image $\operatorname{Im}(\alpha')$ contains
$\alpha'((0,0))=[-4H]=[-D]$, since $H=DF^2$.
Multiplication by $[-D]$ therefore gives a bijection
between its positive and negative classes.
By Lemma~\ref{image-torsors} applied to $E'$, every positive
class has a square-free representative either $\delta$
or $2\delta$, where $\delta>0$ and $\delta\mid D$.
Let $J_0$ and $J_1$ denote the corresponding sets of
positive classes. Then
$|\operatorname{Im}(\alpha')|
=2\bigl(|J_0|+|J_1|\bigr)$.

For a class $[\delta]\in J_0$, put
$\mathbf x=\varphi(\delta)$ as above.
The associated quartic is
$N^2=\delta u^4-4\frac{H}{\delta}v^4$,
with $\gcd(u,v)=1$.
Applying Lemma~\ref{lem-odd-exponent-descent-conditions}
with $\lambda=1$ and $\mu=-4$, we obtain
$A\mathbf x=S\mathbf x$,
or equivalently $Q\mathbf x=\mathbf0_n$.
Thus $|J_0|\le|\ker Q|=2^{k-\eta}$.

For a class $[2\delta]\in J_1$, put
$\mathbf x=\varphi(\delta)$.
The associated quartic is
$N^2=2\delta u^4-2\frac{H}{\delta}v^4$,
with $\gcd(u,v)=1$.
Applying Lemma~\ref{lem-odd-exponent-descent-conditions}
with $\lambda=2$ and $\mu=-2$, we obtain
$A\mathbf x=\mathbf b+S\mathbf x$,
or equivalently $Q\mathbf x=\mathbf b$.
Its solution set is either empty or a coset of $\ker Q$.
Consequently, $|J_1|\le2^{k-\eta}$.

Suppose first that $H\equiv3$ or $15\pmod{16}$.
Then $\eta=1$, and hence
$|\operatorname{Im}(\alpha')| \le2\bigl(2^{k-1}+2^{k-1}\bigr)
=2^{k+1}.$ Now suppose that $H\equiv5$ or $13\pmod{16}$. Here $\eta=0$ and $H\equiv5\pmod8$. We claim that $J_1$ is empty.
Indeed, a primitive solution of $N^2=2\delta u^4-2\frac H\delta v^4$ has $N$ even.
Reduction modulo $4$ shows that $u$ and $v$ have the same
parity, so primitivity forces both to be odd.
Since $\delta$ is odd, $\frac H\delta\equiv H\delta\equiv5\delta\pmod8.$
Therefore $\delta-\frac H\delta\equiv4\pmod8.$
Reducing the quartic modulo $16$ gives $N^2\equiv
2\left(\delta-\frac H\delta\right)
\equiv8\pmod{16},$ which is impossible for a square.
Thus $J_1=\varnothing$, and
$|\operatorname{Im}(\alpha')|
=2|J_0|\le2^{k+1}.$
Writing $r=\operatorname{rank}E_H(\Q)$,
Proposition~\ref{rank proposition} now gives
$$2^r = \frac{|\operatorname{Im}(\alpha)|\,  |\operatorname{Im}(\alpha')|}{4} \le \frac{2^k\,2^{k+1}}4 =2^{2k-1}.$$
Consequently, $\operatorname{rank}E_H(\Q)\le2k-1.$

Finally, if $\operatorname{rank}_{\F_2}A=n-1$,
then $k=1$ and $\operatorname{rank}E_H(\Q)\le1$.
If the Mordell--Weil parity conjecture also holds for $E_H$,
then $(-1)^{\operatorname{rank}E_H(\Q)} =\omega(E_H)=-1$.
The rank is therefore odd and hence equals one.

\item Let $E=E_{-H}: y^2=x^3-Hx, ~ E': y^2=x^3+4Hx,$
and let $\alpha,~\alpha'$ be the corresponding descent maps.
Put $Q=A+S, ~
\eta=\begin{cases}
0,&H\equiv1\pmod4,\\
1,&H\equiv3\pmod4.
\end{cases}$
Since $H\equiv D\pmod4$, Lemma~\ref{lem-odd-D-matrix-rank}
gives $\dim\ker A=k, ~ \dim\ker Q=k-\eta.$

We first prove that $|\operatorname{Im}(\alpha)|
\le2^{k+1-\eta}$.
The image $\operatorname{Im}(\alpha)$ contains
$\alpha((0,0))=[-H]=[-D]$.
Multiplication by $[-D]$ therefore gives a bijection
between its positive and negative classes.
By Lemma~\ref{image-torsors}, every positive class has
a positive square-free representative $\delta\mid D$
for which $N^2=\delta u^4-\frac{H}{\delta}v^4$
has an integer solution satisfying $\gcd(u,v)=1$.
Write $\mathbf x=\varphi(\delta)=(x_1,\ldots,x_n)^{\mathsf T}$,
using \eqref{varphi-def} with $L=D$ and the
prime order in \eqref{matrix A_D}. Applying Lemma~\ref{lem-odd-exponent-descent-conditions}
with $\lambda=1$ and $\mu=-1$, we obtain
$A\mathbf x=S\mathbf x$,
or equivalently $Q\mathbf x=\mathbf0_n$.
Thus the positive classes inject into $\ker Q$,
whose dimension is $k-\eta$.
Consequently, $|\operatorname{Im}(\alpha)|
\le2|\ker Q|=2^{k+1-\eta}$.

We next bound $|\operatorname{Im}(\alpha')|$.
By Lemma~\ref{image-torsors} applied to $E'$, every class
has a signed square-free representative supported on
the primes dividing $2D$ and an associated primitive quartic.
A negative representative would give a quartic with both
coefficients negative. Thus every class is positive.
Let $J_0$ and $J_1$ denote the sets of classes represented
by $\delta$ and $2\delta$, respectively, where
$\delta>0$ and $\delta\mid D$.
Then $|\operatorname{Im}(\alpha')|=|J_0|+|J_1|$.

For a class $[\delta]\in J_0$, put
$\mathbf x=\varphi(\delta)$ as above.
The associated quartic is
$N^2=\delta u^4+4\frac{H}{\delta}v^4$,
with $\gcd(u,v)=1$.
Lemma~\ref{lem-odd-exponent-descent-conditions},
with $\lambda=1$ and $\mu=4$, gives
$A\mathbf x=\mathbf0_n$.
Thus $|J_0|\le|\ker A|=2^k$.

For a class $[2\delta]\in J_1$, put
$\mathbf x=\varphi(\delta)$.
The associated quartic is
$N^2=2\delta u^4+2\frac{H}{\delta}v^4$,
with $\gcd(u,v)=1$.
Lemma~\ref{lem-odd-exponent-descent-conditions},
with $\lambda=\mu=2$, gives
$A\mathbf x=\mathbf b$.
This system has either no solutions or a coset of
$\ker A$ as its solution set.
Consequently, $|J_1|\le2^k$, and hence
$|\operatorname{Im}(\alpha')|\le2^{k+1}$.

Suppose first that $H\equiv7$ or $15\pmod{16}$.
Then $\eta=1$, so the preceding bounds give
$|\operatorname{Im}(\alpha)|\le2^k$ and
$|\operatorname{Im}(\alpha')|\le2^{k+1}$.
Now suppose that $H\equiv5$ or $9\pmod{16}$.
Here $\eta=0$.
We claim that $J_1$ is empty.
Indeed, suppose that $N^2=2\delta u^4+2\frac H\delta v^4,~
\gcd(u,v)=1.$
The right-hand side is even, so $N$ is even.
Reduction modulo $4$ and primitivity force $u$ and $v$
to be odd. Since odd fourth powers are congruent to $1$ modulo $16$, reduction modulo $32$ gives $N^2\equiv
2\left(\delta+\frac H\delta\right)\pmod{32}.$
An even square is congruent to $0$, $4$, or $16$ modulo $32$.
Hence $\delta+\frac H\delta
\equiv0,2,\text{ or }8\pmod{16}.$
For odd integers whose product is $H$, a direct calculation
modulo $16$ gives $\delta+\frac H\delta\equiv0\text{ or }8\pmod{16}$ which implies  $H\equiv7\text{ or }15\pmod{16},$
whereas $\delta+\frac H\delta\equiv2\pmod{16}$ implies $H\equiv1\text{ or }13\pmod{16}.$ 
These possibilities contradict $H\equiv5$ or $9\pmod{16}$.
Thus $J_1=\varnothing$, and $|\operatorname{Im}(\alpha)|\le2^{k+1}, ~|\operatorname{Im}(\alpha')|\le2^k.$
Writing $r=\operatorname{rank}E_{-H}(\Q)$,
Proposition~\ref{rank proposition} gives, in either case,
$$2^r = \frac{|\operatorname{Im}(\alpha)|\,
      |\operatorname{Im}(\alpha')|}{4} \le \frac{2^{2k+1}}4 =2^{2k-1}.$$
Therefore $\operatorname{rank}E_{-H}(\Q)\le2k-1.$

Finally, if $\operatorname{rank}_{\F_2}A=n-1$,
then $k=1$ and $\operatorname{rank}E_{-H}(\Q)\le1$.
If the Mordell--Weil parity conjecture also holds for $E_{-H}$, then $(-1)^{\operatorname{rank}E_{-H}(\Q)}
=\omega(E_{-H})=-1$. The rank is therefore odd and hence equals one.

\item Let $E=E_{-2H}:y^2=x^3-2Hx$ and
$E':y^2=x^3+8Hx$, and let $\alpha,\alpha'$
be the corresponding descent maps.
Put $C=A+B$ and $T=A+B+S$.
Then $\mathcal R_D=[T\mid\mathbf s]$ by
\eqref{eq:augmented-RD}.
By Lemma~\ref{lem:matrix-identities},
$A+A^{\mathsf T}=S+\mathbf s\mathbf s^{\mathsf T}$.
Since $B$ is symmetric,
$T=C^{\mathsf T}+\mathbf s\mathbf s^{\mathsf T}$.
Applying Lemma~\ref{lem-augmented-matrix-rank}
with $X=C^{\mathsf T}$ and $\mathbf u=\mathbf v=\mathbf s$,
we obtain
$\operatorname{rank}_{\F_2}
\binom{C}{\mathbf s^{\mathsf T}}
=\operatorname{rank}_{\F_2}\mathcal R_D=\rho_D$.
Consequently, the rank--nullity theorem gives
\begin{equation}
\label{eq-odd-exponent-kernel-dimensions}
\dim\ker\mathcal R_D=\kappa+1,\quad
\dim\ker
\begin{pmatrix}
C\\ \mathbf s^{\mathsf T} \end{pmatrix}
=\kappa.
\end{equation}

We first prove that $|\operatorname{Im}(\alpha)|
\le2^{\kappa+2}$.
By Lemma~\ref{image-torsors}, every class in this image
has a signed square-free representative supported on
the primes dividing $2D$.
The image contains $\alpha((0,0))=[-2H]=[-2D]$.
Multiplication by $[-2D]$ gives a bijection between
the classes with odd and even square-free representatives.

Let $I_0$ denote the set of classes with odd square-free
representatives. Then
$|\operatorname{Im}(\alpha)|=2|I_0|$.
Every class in $I_0$ has a unique representative
$d=(-1)^t\delta$, where $t\in\{0,1\}$ and
$\delta>0$ with $\delta\mid D$.
Write $\mathbf x=\varphi(\delta)=(x_1,\ldots,x_n)^{\mathsf T}$,
using \eqref{varphi-def} with $L=D$ and the
prime order in \eqref{matrix A_D}. The associated quartic is
$N^2=(-1)^t\delta u^4+(-1)^{t+1}2\frac{H}{\delta}v^4$,
with $\gcd(u,v)=1$.
Applying Lemma~\ref{lem-odd-exponent-descent-conditions}
with $\lambda=(-1)^t$ and $\mu=(-1)^{t+1}2$, we obtain
\[
(A\mathbf x)_i=
\begin{cases}
b_i+(t+1)s_i,&x_i=1,\\
ts_i,&x_i=0.
\end{cases}
\]
Equivalently, $T\mathbf x+t\mathbf s=\mathbf0_n$.
Thus $[(-1)^t\delta]\mapsto\binom{\varphi(\delta)}{t}$,
where $t$ is read in $\F_2$,
injects $I_0$ into $\ker\mathcal R_D$.
By \eqref{eq-odd-exponent-kernel-dimensions},
$|I_0|\le2^{\kappa+1}$.
Consequently,
$|\operatorname{Im}(\alpha)|=2|I_0|\le2^{\kappa+2}$.

We next prove that $|\operatorname{Im}(\alpha')|
\le2^{\kappa+1}$.
By Lemma~\ref{image-torsors} applied to $E'$, every class
has a signed square-free representative supported on
the primes dividing $2D$ and an associated primitive quartic.
A negative representative would give a quartic with both
coefficients negative. Thus every class is positive. The image contains $\alpha'((0,0))=[8H]=[2D]$.
Multiplication by $[2D]$ therefore gives a bijection
between its odd and even classes.
Let $J_0$ denote the set of classes with odd square-free
representatives. Then
$|\operatorname{Im}(\alpha')|=2|J_0|$.

For $[\delta]\in J_0$, where $\delta>0$ and $\delta\mid D$,
put $\mathbf x=\varphi(\delta)$ as above.
The associated quartic is
$N^2=\delta u^4+8\frac{H}{\delta}v^4$,
with $\gcd(u,v)=1$.
Lemma~\ref{lem-odd-exponent-descent-conditions},
with $\lambda=1$ and $\mu=8$, gives
$A\mathbf x=B\mathbf x$,
or equivalently $C\mathbf x=\mathbf0_n$.
If $u$ were even, the condition $\gcd(u,v)=1$
would force $v$ to be odd. Since $\delta$ and
$H/\delta$ are odd, reduction modulo $16$ would give
$N^2=\delta u^4+8\frac{H}{\delta}v^4\equiv8\pmod{16}.$
This is impossible, since a square is congruent to
$0$, $1$, $4$, or $9$ modulo $16$.
Hence $u$ is odd, and the quartic equation shows
that $N$ is odd. Reducing the equation modulo $8$
now gives $\delta\equiv1\pmod8$, and consequently
$\delta\equiv1\pmod4$.
Since $s_i=1$ precisely when $\ell_i\equiv3\pmod4$,
the factorisation $\delta=\prod_{i=1}^n\ell_i^{x_i}$
shows that $\sum_{i=1}^n s_ix_i$ is even.
Thus $\mathbf s^{\mathsf T}\mathbf x=0$
in $\F_2$.
Consequently, $[\delta]\mapsto\varphi(\delta)$ injects
$J_0$ into $\ker\binom{C}{\mathbf s^{\mathsf T}}$.
By \eqref{eq-odd-exponent-kernel-dimensions},
this kernel has dimension $\kappa$.
Therefore $|J_0|\le2^\kappa$, and
$|\operatorname{Im}(\alpha')|=2|J_0|\le2^{\kappa+1}$.
Writing $r=\operatorname{rank}E_{-2H}(\Q)$,
Proposition~\ref{rank proposition} gives
$2^r=|\operatorname{Im}(\alpha)|
|\operatorname{Im}(\alpha')|/4$.
The preceding bounds yield
$2^r\le2^{\kappa+2}2^{\kappa+1}/4=2^{2\kappa+1}$.
Consequently,
$\operatorname{rank}E_{-2H}(\Q)
\le2\kappa+1=2(n-\rho_D)+1$.

If $\rho_D=n$, then $\kappa=0$ and $r\le1$.
If the Mordell--Weil parity conjecture also holds for $E$,
then $(-1)^r=\omega(E)=-1$.
The rank is therefore odd and hence equals one.

\end{enumerate}
\end{proof}

\section{Numerical Examples}\label{Numerical Examples}

Write $E_a:y^2=x^3+ax$.
For a prime support $D=\ell_1\cdots\ell_n$, put $k=n-\operatorname{rank}_{\F_2}\mathcal A_D,~
\rho_D=\operatorname{rank}_{\F_2}\mathcal R_D,~
\kappa=n-\rho_D.$
% All matrix ranks are taken over $\F_2$.
% Every curve listed below has global root number $-1$.
The tables illustrate the unconditional upper bounds
provided by the corresponding theorems.

In Table~\ref{tab:thm-negative-even}, write
$r_D=\operatorname{rank}E_{-2D}(\Q)$. The column headed \emph{Under parity} gives a set containing
$r_D$ under the Mordell--Weil parity conjecture. The square witnesses in the last column establish rank one
unconditionally by Theorem~\ref{Thm negative 2D}\textup{(3)}. The computations in this section were performed using SageMath~\cite{sage}.

\begin{table}[!htbp]
\centering
\setlength{\tabcolsep}{3pt}
\scriptsize

\begin{minipage}[t]{0.49\textwidth}
\centering
\begin{tabular}{@{}lrrrrc@{}}
\toprule
$D$ & $D\bmod16$ & $n$ &
$\operatorname{rank}\mathcal A_D$ & $k$ &
\shortstack{Upper bound\\$2k-1$} \\
\midrule
$5\cdot7=35$         & 3  & 2 & 1 & 1 & 1 \\
$3\cdot7=21$         & 5  & 2 & 1 & 1 & 1 \\
$3\cdot31=93$        & 13 & 2 & 1 & 1 & 1 \\
$3\cdot7\cdot19=399$ & 15 & 3 & 1 & 2 & 3 \\
\bottomrule
\end{tabular}
\caption{Rank bounds for $E_D$
(Theorem~\ref{Thm Positive D}).}
\label{tab:thm-positive-odd}
\end{minipage}
\hfill
\begin{minipage}[t]{0.49\textwidth}
\centering
\begin{tabular}{@{}lrrrrc@{}}
\toprule
$D$ & $D\bmod16$ & $n$ &
$\operatorname{rank}\mathcal A_D$ & $k$ &
\shortstack{Upper bound\\$2k-1$} \\
\midrule
$3\cdot7=21$         & 5  & 2 & 1 & 1 & 1 \\
$7\cdot17=119$       & 7  & 2 & 1 & 1 & 1 \\
$3\cdot7\cdot37=777$ & 9  & 3 & 1 & 2 & 3 \\
$3\cdot5=15$         & 15 & 2 & 1 & 1 & 1 \\
\bottomrule
\end{tabular}
\caption{Rank bounds for $E_{-D}$
(Theorem~\ref{Thm Negative D}).}
\label{tab:thm-negative-odd}
\end{minipage}

\end{table}

% \begin{table}[!htbp]
%   \centering
%   \small
%   \begin{tabular}{@{}lrrrrc@{}}
%     \toprule
%     $D$ & $D\bmod16$ & $n$ &
%     $\operatorname{rank}\mathcal A_D$ & $k$ &
%     \shortstack{Upper bound\\$2k-1$} \\
%     \midrule
%     $5\cdot7=35$         & 3  & 2 & 1 & 1 & 1 \\
%     $3\cdot7=21$         & 5  & 2 & 1 & 1 & 1 \\
%     $3\cdot31=93$        & 13 & 2 & 1 & 1 & 1 \\
%     $3\cdot7\cdot19=399$ & 15 & 3 & 1 & 2 & 3 \\
%     \bottomrule
%   \end{tabular}
%   \caption{Rank bounds for $E_D$
% (Theorem~\ref{Thm Positive D}).}
%   \label{tab:thm-positive-odd}
% \end{table}

% \begin{table}[!htbp]
%   \centering
%   \small
%   \begin{tabular}{@{}lrrrrc@{}}
%     \toprule
%     $D$ & $D\bmod16$ & $n$ &
%     $\operatorname{rank}\mathcal A_D$ & $k$ &
%     \shortstack{Upper bound\\$2k-1$} \\
%     \midrule
%     $3\cdot7=21$         & 5  & 2 & 1 & 1 & 1 \\
%     $7\cdot17=119$       & 7  & 2 & 1 & 1 & 1 \\
%     $3\cdot7\cdot37=777$ & 9  & 3 & 1 & 2 & 3 \\
%     $3\cdot5=15$         & 15 & 2 & 1 & 1 & 1 \\
%     \bottomrule
%   \end{tabular}
%   \caption{Rank bounds for $E_{-D}$
% (Theorem~\ref{Thm Negative D}).}
%   \label{tab:thm-negative-odd}
% \end{table}
\vspace{-0.8cm}
\begin{table}[!htbp]
  \centering
  \small
  \setlength{\tabcolsep}{4pt}
  \begin{tabular}{@{}lrrrccc@{}}
    \toprule
    $D$ & $n$ & $\rho_D$ & $\kappa$ &
    \shortstack{Upper bound $2\kappa+1$} &
    \shortstack{Under parity $r_D\in$} &
    $(\delta,\,\varepsilon,\,h,\,m)$ \\
    \midrule
    $3\cdot11=33$        & 2 & 2 & 0 & 1 & $\{1\}$
      & $(33,\,1,\,2,\,1)$ \\
    $5\cdot17=85$        & 2 & 2 & 0 & 1 & $\{1\}$
      & $(1,\,-1,\,1,\,13)$ \\
    $3\cdot5\cdot7=105$  & 3 & 3 & 0 & 1 & $\{1\}$
      & $(15,\,1,\,1,\,1)$ \\
    $3\cdot5\cdot29=435$ & 3 & 2 & 1 & 3 & $\{1,3\}$
      & \textemdash \\
    \bottomrule
  \end{tabular}
  \caption{Rank bounds for $E_{-2D}$
(Theorem~\ref{Thm negative 2D}).}
  \label{tab:thm-negative-even}
\end{table}

\vspace{-0.95cm}
\begin{table}[!htbp]
  \centering
  \small
  \setlength{\tabcolsep}{4pt}
  \begin{tabular}{@{}lrrlcc@{}}
    \toprule
    $H$ & $D$ & $H\bmod16$ & Curve & Matrix rank &
    \shortstack{Upper bound} \\
    \midrule
    $3^3\cdot29=783$           & 87   & 15 & $E_{783}$
      & $\operatorname{rank}\mathcal A_D=1$ & 1 \\
    % $3^3\cdot13\cdot43=15093$ & 1677 & 5  & $E_{15093}$
    %   & $\operatorname{rank}\mathcal A_D=1$ & 3 \\
    $3^3\cdot11=297$          & 33   & 9  & $E_{-297}$
      & $\operatorname{rank}\mathcal A_D=1$ & 1 \\
    % $3^3\cdot5^3=3375$        & 15   & 15 & $E_{-3375}$
    %   & $\operatorname{rank}\mathcal A_D=1$ & 1 \\
    $3^3\cdot5\cdot7=945$     & 105  & 1  & $E_{-1890}$
      & $\rho_D=3$ & 1 \\
    $3^3\cdot5\cdot67=9045$   & 1005 & 5  & $E_{-18090}$
      & $\rho_D=2$ & 3 \\
    \bottomrule
  \end{tabular}
  \caption{Examples for
Theorem~\ref{thm odd prime exponent extension}.}
  \label{Table for Theorem D}
\end{table}

\bibliographystyle{alpha}
\bibliography{references}

\end{document}